\documentclass{article}
\usepackage{array}
\usepackage{amsmath}
\usepackage{amssymb}
\usepackage{amsthm}
\usepackage{amsfonts}
\usepackage{graphicx}
\usepackage{jca}
\usepackage[toc,page]{appendix}
\usepackage{maplestd2e}
\theoremstyle{definition}

\DefineParaStyle{Maple Heading 1}
\DefineParaStyle{Maple Text Output}
\DefineParaStyle{Maple Dash Item}
\DefineParaStyle{Maple Bullet Item}
\DefineParaStyle{Maple Normal}
\DefineParaStyle{Maple Heading 4}
\DefineParaStyle{Maple Heading 3}
\DefineParaStyle{Maple Heading 2}
\DefineParaStyle{Maple Warning}
\DefineParaStyle{Maple Title}
\DefineParaStyle{Maple Error}
\DefineCharStyle{Maple Hyperlink}
\DefineCharStyle{Maple 2D Math}
\DefineCharStyle{Maple Maple Input}
\DefineCharStyle{Maple 2D Output}
\DefineCharStyle{Maple 2D Input}

\usepackage[top=.9in,bottom=.9in]{geometry}
 \numberwithin{equation}{section}
 
\DefineParaStyle{Maple Output}
\DefineCharStyle{2D Comment}
\DefineCharStyle{2D Math}
\DefineCharStyle{2D Output}
\begin{document}
\let\labeldefs\iftrue
\let\labeldefs\iffalse
\pagestyle{plain}
\numberwithin{equation}{section}
\begin{flushleft} 
\vskip 0.3 in  
\centerline{\bf Variations on a Hypergeometric Theme} \vskip .3 in
\vskip .2in
\centerline{ Michael Milgram\footnote{mike@geometrics-unlimited.com}}

\centerline{Consulting Physicist, Geometrics Unlimited, Ltd.}
\centerline{Box 1484, Deep River, Ont. Canada. K0J 1P0}
\centerline{}

\centerline{}
\centerline{revised Oct. 11, 2018: Appendix C revised and corrected}
\vspace{.5cm}
\centerline{}
\centerline{}
\vskip .3in
\centerline{\bf Abstract}\vskip .1in

The question was asked: Is it possible to express the function

\begin{equation} \tag{1.1}
h(a)\equiv\,{_4F_3}(a,a,a,a;2a,a+1,a+1;1)
\label{question}
\end{equation}
\labeldefs{\begin{flushright}question\end{flushright}}\fi

in closed form \cite{Mazac}? After considerable analysis, the answer appears to be ``no", but during the attempt to answer this question, a number of interesting (and unexpected) related results were obtained, either as specialized transformations, or as closed-form expressions for several related functions. The purpose of this paper is to record and review both the methods attempted and the related identities obtained, the former for their educational merit, the latter because they do not appear to exist in the literature. Specifically, new $_4F_3(1)$, $_5F_6(1)$ and generalized Euler sums (those containing digamma functions) are presented along with a detailed discussion of the methods used to obtain them. \newline

\section{Introduction and Motivation} 

Always keeping the question posed by \eqref{question} in mind, this paper is an exposition, or review, of well-established techniques which are not always employed in the modern literature. Thus, relationships that have been derived by long, perhaps convoluted developments, often found in the literature, turn out to be special cases of well-established, but rarely-used identities. For example Miller and Paris \cite[Eq. (1.7)]{Miller&Paris} obtain a representation of $_3F_2(a,b,n;c,m;1)$ in terms of two finite sums (see \eqref{Eq1p7a}) that is easily shown to be a special case of a known 3-part transformation among various $_3F_2(1)$ (see \eqref{Eq1p7b}).\newline 

The approach explored here revolves about the use of two, three and four-part transformations between and among general hypergeometric functions of argument unity. In the first two appendices, the most general transformations among $_{p+1}F_{p}(1) (p\le3)$ functions that appear to have potential to be profitably applied to \eqref{question} are developed. A third Appendix revisits the long-ago-established, but rarely invoked, notation for 3-part transformations among $_3F_2(1)$. A final Appendix collects a number of useful results that are either scattered throughout the literature, or required as lemmas for the derivations contained throughout the text. In the main text, (Sections \ref{sec:Consult}, \ref{sec:Xform} and \ref{sec:Evol}) each general transformation identified in an appendix is applied to the problem \eqref{question}, in a corresponding subsection. Thus it is recommended that each appendix be consulted along with its corresponding section to obtain a clear picture - that is, {\bf read (at least) the (first two) appendices first!} In general each such application yields a new representation for $h(a)$, or, in a few cases, a functional relationship between $h(a)$ and $h(1-a)$. An observant reader will note that many of the results are symmetric under the exchange $a\rightarrow 1-a$; I have attempted to present each representation is such a way that this symmetry is clear.\newline

Having failed to solve \eqref{question} in the first few sections, a reduced problem is considered in Section \ref{sec:reduce}.  When that in turn fails to produce an answer, a novel technique of dubious mathematical rigour is attempted in Section \ref{sec:hail} with partial success. In the final Sections \ref{sec:Consq} and \ref{sec:new}, the various representations that were derived as described above are compared to one another, which has the unexpected effect of yielding some new representations, including what may be several new $_4F_3(1)$ and sums involving digamma functions. The discovery of one such, (i.e. \eqref{4F3}) in turn inspired a detailed study of related forms, which in their own turn, resulted in more new identities and closed forms for a variety of $_4F_3(1)$, $_5F_6(1)$ and related sums.\newline

Throughout, all variables are assumed to be complex ($\mathfrak{C}$) except in special instances as noted, with the overall exception that the symbols $k$, $m$ and $n$, are positive integers ($k,m,n\in\mathbb{N}$). $\Psi(a), \Psi^{\prime}(a)$ and $\Psi^{(n)}(a)$ are respectively the digamma function, its first and $n^{th}$ derivative, $n>1$. As usual, $\gamma$ represents the Euler constant.  In many cases, the use of a limit procedure is indicated by the symbol ``$\rightarrow$", whereas symbolic replacement (redefinition) is indicated by the use of  ``$:=$". Although the notation used throughout is the typographically convenient $_pF_q(\lbrace a \rbrace;\lbrace b \rbrace;1)$ rather than the less convenient $_pF_q(\substack{\lbrace\text{a}\rbrace\\\lbrace\text{b}\rbrace}|1)$, textual references ``top" and ``bottom" are applied to the parameters. The meaning should be clear. It is tacitly required that the combination of parameters appearing in any of the hypergeometric identities is such that the parametric excess (real part of (sum of bottom minus sum of top) parameters) is positive. Finally, it is noted that the methods outlined here, although applied to a simple hypergeometric function of one parameter (variable), should have application to similar problems involving more parameters.

\section{Method 1: Consult both the literature and computer algorithms}\label{sec:Consult}

Reference literature on the subject of $_4F_3(1)$ is very limited.  Although sums closely related to \eqref{question} appear in the standard reference works \cite[Eqs. 7.6.2(13) and 7.10.2(6)]{prudnikov},\cite[Chapter 2]{Slater} and \cite[Eq. (10.39.4) with $a=b$]{Hansen} the exact sum sought does not appear. Similarly, computer codes (e.g. Maple, Mathematica) incorporate  telescoping methods (e.g. \cite{Gosper}) which can be used to attempt to simplify any desired hypergeometric function; this approach did not yield a closed form for \eqref{question}. In retrospect, this is unsurprising, since all of the representations in this case involve digamma functions whose underlying structure is not strictly hypergeometric in nature - a requirement of the telescoping algorithms.

\section{Method 2: Consider a $_4F_3(1)$ Transformation} \label{sec:Xform}

It is well known that the much utilized Watson/Whipple/Dixon identities for a $_3F_2(1)$ can be obtained from one another using any of the ten independent 2-part Thomae relations \cite{Bailey1935}, so knowledge of one will yield another. Perhaps the query function \eqref{question} can be related to other (known) $_4F_3(1$) by a simple transformation?\newline 

Explicit two-part transformations between $_4F_3(1)$ \`a la Thomae are inapplicable (there are perhaps 160 when the $_4F_3(1)$ is Saalsch\"utzian with integral parameters \cite{Raynal}). However, Miller \cite{Miller} has obtained two potentially relevant 4-part transformations applicable to $_4F_3(1)$ with arbitrary parameters - see Appendix \eqref{Miller1} and \eqref{Miller2}. Further, in Appendix A, a new 4-part transformation of the same form is  obtained from the analytic continuation of Meijer's G-function with respect to its argument $z\rightarrow 1/z$ when $z=1$. See \eqref{GF_Imag}.\newline

Other (inapplicable) possibilities include:
\begin{itemize}
\item an explicit result for a special Saalsch\"utzian $_4F_3(1)$ when a top parameter is a negative integer (see \eqref{WhipXf1} and \cite{Choi_Honam}) - note that two equations in \cite{Choi_Honam} are misprinted: in the right-hand side of Eq.(1.10) the denominator parameter $a$ has the wrong sign, and in Eq. (2.2) replace $c$ by $f$;
\item the generalized Minton-Karlsson reduction \cite{Minton}, \cite{Karlsson}, \cite{Sriv_Minton_K} when a top parameter exceeds a bottom parameter by a positive integer, and;
\item \cite[ Eq. (20)]{Gottschalk} where Gottschalk and Maslen propose a reduction formula for a general $_pF_q$ that outwardly appears to apply to \eqref{question} except that it is valid only for $q>p$; for the case $p=q+1$ the arbitrary parameter $b$ they introduce (potentially identified as $2a$ in \eqref{question}) appears as a top - not bottom - parameter (see Appendix A.2).\end{itemize}

With this background, it is instructive to carefully investigate the first of Miller's transformations \eqref{Miller1} as reproduced in Appendix A. The following detailed sequence, using computer algebra, should be studied as a template for all similar calculations herein. It is equivalent to older, very complicated general results claimed in the literature (for example, see \cite{SaiSax} and references therein; however Eq. (2.1) in that article does not satisfy numerical tests). In \eqref{Miller1}, let $e=2a,\,f=a+1$ and $g=a+1$ to obtain

\begin{align} \nonumber
_4F_3&(a,b,c,d;\,2\,a,1+a,1+a;\,1)=  \Gamma \left( 1-d \right) \Gamma \left( 2\,a+1 \right) \Gamma \left( 1+a \right) \\ \nonumber
& \hspace{-.3cm} \times \left( {\frac {{_4F_3(b,b-a,b-a,1+b-2\,a;\,1+b-d,1+b-a,1+b-c;\,1)} }{2\,\Gamma \left(1+ a-b \right) 
\Gamma \left( c \right) \Gamma \left( 2\,a-b \right)  \left( a-b \right) \Gamma \left( 1+b-d \right) }\Gamma \left( c-b \right)} \right. \\ \nonumber
& \left. +{\frac { {_4F_3(c,c-a,c-a,1+c-2\,a;\,1+c-d,1+c-a,1+c-b;\,1)}}{2\,\Gamma \left( a-c+1 \right)  \left( a-c \right) \Gamma \left( 2\,a-c \right) \Gamma \left( 1+c-d \right) 
\mbox{}\Gamma \left( b \right) }\Gamma \left( b-c \right)}
\right. \\ 
& \left. +\,{\frac {\Gamma \left( b-a \right) \Gamma \left( c-a \right) }{2\,\Gamma \left( c \right) \Gamma \left( 1+a-d \right) \Gamma \left( b \right) }} \right)
\label{MStep1}
\end{align}

With $b=d$ one finds
\begin{align} \nonumber
\displaystyle \mbox{$_4$F$_3$}&(a,c,d,d;\,2\,a,1+a,1+a;\,1)=\frac {\Gamma \left( 1+a \right) \Gamma \left( 2\,a+1 \right) \Gamma \left( 1-d \right) }{\Gamma \left( c \right) \Gamma \left( d \right) } \\ \nonumber
&\left( {\frac {\sin \left( \pi\, \left( c-a \right)  \right) \sin \left( \pi\, \left( 2\,a-c \right)  \right) }{2\,\pi\,\sin \left( \pi\, \left( d-c \right)  \right) }\sum _{k=0}^{\infty }{\frac {\Gamma \left( 1+c-2\,a+k \right) \Gamma \left( c-a+k \right) \Gamma \left( c+k \right) }{ \left( a-c-k \right)  \Gamma \left( 1+c-d+k \right) ^{2}
\mbox{}\Gamma \left( k+1 \right) }}} \right. \\ \nonumber
&\left. +{\frac {\sin \left( \pi\, \left( a-d \right)  \right) \sin \left( \pi\, \left( -d+2\,a \right)  \right) }{2\,\pi\,\sin \left( \pi\, \left( d-c \right)  \right) }\sum _{k=0}^{\infty }{\frac {\Gamma \left( 1+d-2\,a+k \right) \Gamma \left( d-a+k \right) \Gamma \left( k+d \right) }{ \left( a-d-k \right)  \Gamma \left( k+1 \right)^{2}
\mbox{}\Gamma \left( k+1-c+d \right) }}} \right. \\ 
&\left.+{\frac {\Gamma \left( d-a \right) \Gamma \left( c-a \right) }{2\,\Gamma \left( 1+a-d \right) }} \right) 
\label{MStep2}
\end{align}
and further, after taking the limit $c\rightarrow d$
\begin{flalign} \label{Mstep3} \nonumber 
\frac {2\,{\pi}^{2}\Gamma \left( d \right) }{\Gamma \left( 1+a \right) \Gamma \left( 1-d \right) \Gamma \left( 2\,a+1 \right) 
} &  {  {_4F_3(a,d,d,d;\,2\,a,1+a,1+a;\,1)}} =   -\pi\,\sin \left( 3\,\pi\,a-2\,\pi\,d \right) \sum _{k=0}^{\infty }S \left( k \right) \\ \nonumber
&\hspace{-5cm} + \left( 2\,\sum _{k=0}^{\infty }{\frac {S \left( k \right) }{a-d-k}}+\sum _{k=0}^{\infty }S \left( k \right) \Psi \left( 1+d-2\,a+k \right) +\sum _{k=0}^{\infty }S \left( k \right) \Psi \left( 1+d-a+k \right) 
\mbox{} \right. \\ & \hspace{-5cm} \nonumber
\left.
 +\sum _{k=0}^{\infty }\Psi \left( k+d \right) S \left( k \right) -3\,\sum _{k=0}^{\infty }\Psi \left( k+1 \right) S \left( k \right)  \right) \sin \left( \pi\, \left( -d+2\,a \right)  \right) \sin \left( \pi\, \left( a-d \right)  \right) \\ \left.
+\frac {{\pi}^{4}}{\sin^{2} \left( \pi\, \left( a-d \right)  \right) \Gamma \left( 1+a-d \right)^{3}} \right. 
\end{flalign}
where
\begin{equation}
\mapleinline{inert}{2d}{S(k) = GAMMA(k+d)*GAMMA(d-a+k)*GAMMA(1+d-2*a+k)/(GAMMA(k+1)^3*(a-d-k))}{\[\displaystyle S \left( k \right) ={\frac {\Gamma \left( k+d \right) \Gamma \left( d-a+k \right) \Gamma \left( 1+d-2\,a+k \right) }{ \Gamma \left( k+1 \right)  ^{3} \left( a-d-k \right) \,.
\mbox{}}}\]}
\label{S(k)}
\end{equation}
Finally, taking the limit $d\rightarrow a$ eventually gives the following apparently new representation after much cancellation and simplification of terms

\begin{align} \label{Miller1Final}
\mapleinline{inert}{2d}{hypergeom([a, a, a, a], [2*a, 1+a, 1+a], 1)/a^2 = -GAMMA(2*a)*a*(-1+a)*Pi*hypergeom([1, 1, 1, 2-a, 1+a], [2, 2, 2, 2], 1)/(GAMMA(a)^2*sin(Pi*a))+2*Pi*(gamma+Psi(a))^2*GAMMA(2*a)/(GAMMA(a)^2*sin(Pi*a))}{\[\displaystyle {\frac {h(a)}{{a}^{2}}}=a\,\pi\,\left(1-a \right){\frac {\Gamma \left( 2\,a \right)   \,{\mbox{$_5$F$_4$}(1,1,1,1+a,2-a;\,2,2,2,2;\,1)}}{ \Gamma \left( a \right)  ^{2}\sin \left( \pi\,a \right) 
\mbox{}}}+2\,{\frac {\pi\, \left( \gamma+\Psi \left( a \right)  \right) ^{2}\Gamma \left( 2\,a \right) }{ \Gamma \left( a \right)   ^{2}\sin \left( \pi\,a \right) }}\]}\,.
\end{align}

For the case $a=1$ see Section \ref{sec:reduce}; it is also convenient to quote the simple contiguity relation
\begin{align} \nonumber
\displaystyle {\mbox{$_4$F$_3$}(1,1,a+1,1-a;\,2,2,2;\,1)}=a\,{\mbox{$_5$F$_4$}(1,1,1,1+a,1-a;\,2,2,2,2;\,1)} \\
-\left( a-1 \right)\, {\mbox{$_5$F$_4$}(1,1,1,1+a,2-a;\,2,2,2,2;\,1)}\,.
\label{Heq}
\end{align}

With reference to Appendix A, and in a similar fashion, the same series of substitutions and limits applied to Miller's second transformation \eqref{Miller2} surprisingly yields the same result. However, when this sequence of substitutions and limits is similarly applied to the transformation \eqref{GF_Imag}, a different representation in the form of a functional equation emerges:
%
\begin{align} \label{Miller2Final} 
\frac{h(a)}{a^2}  =\frac{h(1-a)}{(1-a)^2}  {\frac {\Gamma \left( 2\,a \right)  \Gamma \left( 1-a \right)^{2}}{  
\Gamma \left( 2-2\,a \right)  \Gamma \left( a \right) ^{2}}}
\mbox{}  
-4\,{\pi}^{2}{\frac {\Gamma \left( 2\,a \right) \cos \left( \pi\,a \right) }
{\sin^2 \!\left( \pi\,a \right)  \Gamma \left( a \right)^{2}} \left( \gamma+\Psi \left( a \right) +\frac{\pi}{2}{ {\cot \left( \pi\,a \right) }} \right) }
\end{align}

\section{Evolution from a $_3F_2(1)$} \label{sec:Evol}

\subsection{The Principle} 
The method employed here is based on the observation that, written as a sum, the hypergeometric function 
\begin{equation}
\mapleinline{inert}{2d}{Hypergeom([a, a, b], [2*a, b+1], 1) = GAMMA(2*a)*b*(Sum(GAMMA(a+k)^2/(GAMMA(k+1)*GAMMA(2*a+k)*(b+k)), k = 0 .. infinity))/GAMMA(a)^2}{\[\displaystyle {_3F_2} \left( a,a,b\,;2\,a,b+1\,;1 \right) =b\,{\frac {\Gamma \left( 2\,a \right)}{  \Gamma \left( a \right)^{2}}\sum _{k=0}^{\infty }{\frac { \Gamma \left( a
\mbox{}+k \right)^{2}}{\Gamma \left( k+1 \right) \Gamma \left( 2\,a+k \right)  \left( b+k \right) }}}\]}
\label{3F2}
\end{equation}
yields the following rule after being acted upon by the consecutive operations $\displaystyle \frac{\partial}{\partial\,b}$ followed by ${\it b=a}\;$

\begin{align} \label{Heq1} 
\displaystyle h(a) =
{_3F_2} \left( a,a,a;2\,a,a+1;1 \right) 
\mbox{}-\frac{{a}^{2}}{2}\,{\frac {{_4F_3} \left( a+1,a+1,a+1,a+1\,;a+2,a+2,2\,a+1;1 \right) }{ \left( a+1 \right) ^{2}}}\,.
\end{align}
In the following I represent this sequence of actions by the operator

\begin{equation} \label{MathL}
\mathfrak{L}(b)\equiv b\rightarrow a \frac{\partial}{\partial\,b}.
\end{equation}

Note that the $_4F_3(1)$ on the right-hand side of \eqref{Heq1} is almost related to $h(a)$ by the operation $a\rightarrow a+1$ -  see \eqref{R12}. This result is particularly interesting because the ${_3F_2(1)}$ that appears in \eqref{Heq1} is a limiting case contiguous to Watson's theorem \cite{ChuW} and therefore can be evaluated in closed form:

\begin{equation}
\mapleinline{inert}{2d}{Hypergeom([a, a, a], [2*a, a+1], 1) = (1/4)*2^(2*a)*a*(-Psi(1, (1/2)*a+1/2)+Psi(1, (1/2)*a))*GAMMA(a+1/2)/(sqrt(Pi)*GAMMA(a))}{\[\displaystyle {_3F_2} \left( a,a,a;2\,a,a+1;1 \right) =a\,{\frac {{2}^{2\,a} \left( -\Psi^{\prime} \left( a/2+1/2 \right) +\Psi^{\prime} \left( a/2 \right)  \right) 
\mbox{}\Gamma \left( a+1/2 \right) }{ 4\,\sqrt{\pi}\;\Gamma \left( a \right) }}\]}\,.
\label{Entry}
\end{equation}

In the following sections and sub-sections, various representations for $_3F_2(a,a,b;2a,b+1;1)$ based on different 3-part transformations of a generic function $_3F_2(a,b,c;e,f;1)$ are subjected to the action of the operator $\mathfrak{L}(b)$, to determine if any useful representations for the function $h(a)$ can be found. The various transformations are developed in  Appendix B; as well, a primer on the general subject of 3-part transformations of $_3F_2(1)$ will be found in Appendix C.
\newline
 
\subsection{ Based on \eqref{Eq1p7b}}

For the particular case $h(a)$ in \eqref{Eq1p7b}, based on a result of Miller and Paris, set $f=c+1,b=a,e=2a$ then operate on both sides with $\displaystyle \frac{\partial}{\partial\,c}$ followed by the limit ${c\rightarrow a}\;$ and after all the various limits and substitutions are completed, obtain
\begin{align} \label{H0a} \nonumber
&\displaystyle {\frac {h(a)}{{a}^{2}}}={\frac {\pi\,\Gamma \left( 2\,a \right) 
\mbox{}}{\sin \left( \pi\,a \right) \Gamma \left( a \right) ^{3}}\sum _{k=0}^{\infty }{\frac {\Gamma \left( k+1 \right) \Psi \left( 1-a+k \right) }{\Gamma \left( 1-a+k \right)  \left( 1-a+k \right)^2 }}}
\\ \nonumber
&+{\frac {\Psi \left( 1-a \right) \pi\,\Gamma \left( 2\,a \right) {_3F_2} \left( 1,1,1-a;2-a,2-a;\,1 \right) 
\mbox{}}{\Gamma \left( 2-a \right)  \left( -1+a \right) \Gamma \left( a \right)^{3}\sin \left( \pi\,a \right) }}
\\
&-\,{\frac {\sin \left( 2\,\pi\,a \right)  \left( 2\gamma+\Psi \left( 1-a \right) +\Psi \left( a \right)  \right) \Gamma \left( 1-a \right) 
\mbox{}\pi\,\Gamma \left( 2\,a \right) }{\Gamma \left( a \right)  \sin^2 \left( \pi\,a \right) }}\,.
\end{align}

The series appearing in \eqref{H0a} is convergent by Gauss' test if $\Re(a)<1$. Additionally, the $_3F_2(1)$ appearing is easily identified because it is contiguous to Whipple's theorem \cite[Entry 13] {Milgram447}, thus
\begin{equation}
\mapleinline{inert}{2d}{hypergeom([1, 1, 1-a], [2-a, 2-a], 1) = (1/2*(Psi(1, 1/2-(1/2)*a)+Psi(1, (1/2)*a)-Pi^2/sin((1/2)*Pi*a)^2))*(-1+a)^2}{\[\displaystyle {\mbox{$_3$F$_2$}(1,1,1-a;\,2-a,2-a;\,1)}=1/2\, \left( \Psi^{\prime} \left(1/2-a/2 \right) +\Psi^{\prime} \left(a/2 \right) -{\frac {{\pi}^{2}}{ \sin^2 \left( \pi\,a/2 \right)  }} \right)  \left( -1+a \right) ^{2}\]}\,,
\label{TT}
\end{equation}

yielding the  following after substitution and simplification
\begin{align} \label{Eq4}
&\displaystyle {\frac {h \left( a \right) }{{a}^{2}}}={\frac {\pi\,\Gamma \left( 2\,a \right) 
\mbox{}}{\sin \left( \pi\,a \right)   \Gamma \left( a \right)  ^{3}}\sum _{k=0}^{\infty }{\frac {\Gamma \left( k+1 \right) \Psi \left( 1-a+k \right) }{\Gamma \left( 1-a+k \right)  \left( 1-a+k \right)^2 }}} \hspace{3cm}a<1\\ \nonumber
&+{\frac {\Gamma \left( 2\,a \right) }{\Gamma \left( a \right)  ^{2}} \left( -4\, \gamma\,{\pi}^{2}\,{\frac {\cos \left( \pi\,a \right)}{  \sin^{2} \left( \pi\,a \right) }}-\Psi \left( 1-a \right) \Psi^{\prime} \left(a/2 \right) +2\,\Psi \left( 1-a \right) \Psi^{\prime} \left( a \right) 
-2\,{\pi}^{2}\,{\frac {\cos \left( \pi\,a \right) \Psi \left( a \right) }{\sin^2 \left( \pi\,a \right) }} \right) }
\end{align}

\subsection{Based on \eqref{Shpot3}}

After applying the set of transformations discussed above to the generalization \eqref{Shpot3}  of the result of Shpot and Srivastava \cite{Shpot&Sriv}, \eqref{Miller1Final} re-emerges.
 
\subsection{Based on \eqref{F32_Imag}}

Similarly, after applying the set of substitutions and operations  $c=a,e=2a, f=b+1$ using the transformation rule \eqref{F32_Imag} as a basis, one finds

\begin{align} \label{FromB7}
&\displaystyle {\mbox{$_3$F$_2$}(a,a,b;\,2\,a,b+1;\,1)}= {\frac {\pi\,\sin \left( 2\,\pi\,a \right) \Gamma \left( b+1 \right) \Gamma \left( 1+b-2\,a \right) \Gamma \left( 2\,a \right) }{ \sin^2 \left( \pi\,a \right)   \Gamma \left( b+1-a \right) ^{2}
 \Gamma \left( a \right) ^{2}}} \\ \nonumber
&+{\frac {{b\,\pi}^{2}\Gamma \left( 2\,a \right) {\mbox{$_3$F$_2$}(1-a,1-a,b+1-2\,a;\,2-2\,a,b+2-2\,a;\,1)}
\mbox{}}{  \sin^2 \left( \pi\,a \right)  \Gamma \left( a \right) ^{4}\Gamma \left( 2-2\,a \right)  \left( b+1-2\,a \right) }}
\end{align}

and, after operating with $\mathfrak{L}(b)$, the final result reduces to \eqref{Miller2Final}.

\subsection{Based on \eqref{Miller2b}}

Again, and following the same procedure as above, use the transformation rule \eqref{Miller2b} as a basis to eventually obtain
\begin{align} \label{FromMiller2b}
&\displaystyle {\frac {h \left( a \right) }{{a}^{2}}}={\frac {4\,{\pi}^{2}\Gamma \left( 2\,a \right)   ^{2} \cos^2 \left( \pi\,a \right) {\mbox{$_3$F$_2$}(1-a,1-a,1-a;\,2-a,2-2\,a;\,1)}
\mbox{}}{ \Gamma \left( a \right)^{4} \left(2\,a  -1\right)  \left( a-1 \right) \sin^2 \left( \pi\,a \right) }}
\\ \nonumber
&-\pi\,a{\frac {  \Gamma \left( 1-a \right) ^{2} \left( a-1 \right) {\mbox{$_4$F$_3$}(1,1,2-a,1+a;\,2,2,2;\,1)}}{\Gamma \left( 1-2\,a \right) \sin \left( \pi\,a \right) }}
\\ \nonumber
& +{\frac {1}{ \Gamma \left( a \right)^{2}} \left( {\frac {\Gamma \left( 2\,a \right) \Gamma \left( 1-a \right)^{2}h \left( 1-a \right) }{\Gamma \left( 2-2\,a \right)  \left( a-1 \right) ^{2}}} 
+{\frac {{\pi}^{4}\cos \left( \pi\,a \right) }{\Gamma \left( 1-2\,a \right)  \sin^4 \left( \pi\,a \right)  }} 
\mbox{} \right) } \,.
\end{align}

The $_3F_2(1)$ appearing in the first term of \eqref{FromMiller2b} can be identified from \eqref{Entry} with $a\rightarrow 1-a$, and when this is done and substituted into \eqref{FromMiller2b} one (ultimately) finds another different but equivalent version of the functional equation \eqref{Miller2Final}:
 
\begin{align} \label{Miller2d} 
&\displaystyle {\frac {\Gamma \left( a \right) h \left( a \right) \sin \left( \pi\,a \right) }{{a}^{2}\Gamma \left( a+1/2 \right) }}=-{\frac {{4}^{2\,a-1}\Gamma \left( a+1/2 \right) \cos \left( \pi\,a \right) h \left( 1-a \right) }{ \left( a-1 \right) ^{2} \left( a-1/2 \right) 
\mbox{}\Gamma \left( a \right) }} 
+{2}^{2\,a} \sqrt{\pi}\cos \left( \pi\,a \right)  \left( \Psi^{\prime} \left(a/2 \right) -2\,\Psi^{\prime} \left( a \right)  \right) 
\mbox{} \\ \nonumber
&-{\pi}^{3/2}a \left( a-1 \right) {2}^{2\,a}{\mbox{$_4$F$_3$}(1,1,2-a,a+1;\,2,2,2;\,1)}\cot \left( \pi\,a \right)
 -{\pi}^{5/2}{2}^{2\,a} \cot^2 \left( \pi\,a \right)  \,.
\end{align}

\subsection{Based on \eqref{Mix1}}

In Appendix B, a curious reduction of a $_3F_2(a,b,c;e,f;1)$ was developed - specifically, an arbitrary parameter $d$ appeared on the right-hand side, when a template $_3F_2(a,b,c;e,f;1)$  was expressed in the form of a 3-part transform. After applying to \eqref{Mix1} the operations discussed above, the following curious representation emerges
\begin{align} \label{Sm1b}
&\displaystyle {\mbox{$_3$F$_2$}(a,a,b;\,2\,a,b+1;\,1)}=-{{2}^{2\,a-1}\,b\,\sqrt{\pi}\,\Gamma \left( a+1/2 \right)\frac {  {\mbox{$_3$F$_2$}(a,a-b,1-a;\,1,1+a-b;\,1)}}{ \left( a-b \right) \Gamma \left( a \right) 
\mbox{}\sin \left( \pi\,a \right) }} \\ \nonumber
&+{b\,\frac {\Gamma \left( 2\,a \right) \sin \left( \pi\, \left(a-d \right)  \right) {\mbox{$_4$F$_3$}(a,a-b,a-b,1-a;\,1,1+a-b,1+a-b;\,1)}}{ \left( a-b \right) ^{2} \Gamma \left( a \right) ^{2}
\mbox{}\sin \left( \pi\,d \right) }} \\ \nonumber
&-{b\,S_{{2}}(a,b)\frac {\Gamma \left( 2\,a \right) \sin \left( \pi\, \left(a -d \right)  \right) \sin \left( \pi\,a \right) }{\pi\,\Gamma \left( a \right)  ^{2}\sin \left( \pi\,d \right) }}
\mbox{}-{\frac {\Gamma \left( 2\,a \right) \Gamma \left( b+1 \right) \Gamma \left( a-b \right) ^{2}\sin \left( \pi\, \left( b-d \right)  \right) }{\Gamma \left( a \right) ^{2}\Gamma \left( 2\,a-b \right) 
\mbox{}\sin \left( \pi\,d \right) }}
\end{align}

where

\begin{align} \label{S2}
\displaystyle { S_{2}}(a,b)\, \equiv & \,\sum _{k=0}^{\infty }{\frac {\Gamma \left( k+1-a \right) \Gamma \left( a+k \right)\left( \Psi \left( k+1-a \right)+\Psi(a+k)-2\Psi(k+1)\right) }{ \left( a-b+k \right) \Gamma \left( k+1 \right)  ^{2}
\mbox{}}}\,.
\end{align}
For the case $d=0$ see \eqref{Smm1H}.\newline


Again, operate with  $ \mathfrak{L}(b)$, to find a very long and complicated representation. After considerable simplification and application of independently obtained results \eqref{4F3} and \eqref{6F5} (see below) as well as \eqref{Miller1Final} and \eqref{Entry} above, the resulting representation is

\begin{align} \label{Sm3d}
\displaystyle {\frac {h \left( a \right) }{{a}^{2}}}=&{\frac {{4}^{a}\Gamma \left( a+1/2 \right) }{6 \sqrt{\pi}\Gamma \left( a \right) }}{S_{1}(a)} \\ \nonumber
&+\frac {  {4}^{a}\Gamma \left( a+1/2 \right) 
\mbox{}}{3\,\Gamma \left( a \right) \sin \left( \pi\,a \right) } 
\left( -{\displaystyle \frac{{\pi}^{5/2}}{6}}+ \left( \gamma+\Psi \left( a \right)  \right) \cot \left( \pi\,a \right) {\pi}^{3/2}+4\, \left( \gamma+\Psi \left( a \right)  \right) ^{2} \sqrt{\pi}
\mbox{} \right)
\end{align}

where the absolutely convergent series $S_{1}(a)$ is defined by (also see \eqref{S1_new})

\begin{align} \label{S1}
\displaystyle { S_{1}(a)}\, \equiv & \,\sum _{k=1}^{\infty }{\frac {\Gamma \left( k+1-a \right) \Gamma \left( a+k \right)\left( \Psi \left( k+1-a \right)+\Psi(a+k)-2\Psi(k)\right) }{k^3\, \Gamma \left( k \right)  ^{2}
\mbox{}}}\,.
\end{align}

Notice that dependence on the parameter $d$ has vanished during the simplification process. Alternatively, from \eqref{Mix1} operate with the following sequence: $f=b+1, e=2a$ followed by the limit $c\rightarrow a$, then $\mathfrak{L}(b)$, and finally $d=a$, to obtain

\begin{align} \label{MdEqa}
{h(a)}=&{{_3F_2}(a,a,a;\,2\,a,a+1;\,1)} +\frac{\pi\, \Gamma \left( 2\,a+1 \right)}{\sin \left( \pi\,a \right)\, \Gamma \left( a \right)   ^{2}}  
 \left( -\frac{1}{2}\,a \left( a-1 \right) {{_4F_3}(1,1,a+1,2-a;\,2,2,2;\,1)} \right. \\ \nonumber
& \left. -\frac{1}{2}\,{a}^{2} \left( a-1 \right)
 {{_5F_4}(1,1,1,a+1,2-a;\,2,2,2,2;\,1)} \right.   
  \left.  + \left( \Psi \left( a \right) +\gamma \right) ^{2}a
+\Psi \left( a \right) +\gamma \right) \,.
\end{align}
When \eqref{Entry} and \eqref{4F3} are applied, \eqref{MdEqa} reduces to \eqref{Miller1Final}.

\subsection{Based on \eqref{Mix2}}
As discussed in the Appendix, the result \eqref{Mix2} also admits an arbitrary parameter $d$ on the right-hand side. After operating on \eqref{Mix2} with the following left-ordered sequence $d=a,b=a, f\rightarrow b+1, e=2\,a$ we find the interim result
\begin{align} \label{Hb3B}
&\displaystyle {\frac {{\mbox{$_3$F$_2$}(a,a,b;\,2\,a,1+b;\,1)}}{b}}=a\,\pi\, \left( a-1 \right){\frac {\Gamma \left( 2\,a \right) \,{\mbox{$_4$F$_3$}(1,2-a,1+a,1+a-b;\,2,2,a-b+2;\,1)}}{ \Gamma \left( a \right) ^{2}
\mbox{}\sin \left( \pi\,a \right)  \left( 1+a-b \right) }} \\ \nonumber
&+{\frac {{2}^{2\,a-1} \sqrt{\pi}\Gamma \left( a+1/2 \right) \Gamma \left( b \right) \Gamma \left( a-b \right) }{\Gamma \left( a \right) \sin \left( \pi\,a \right) \Gamma \left( 1+b-a \right) \Gamma \left( 2\,a-b \right) 
\mbox{}}}- 
{\frac {{2}^{2\,a-1} \sqrt{\pi}\Gamma \left( a+1/2 \right) }{\Gamma \left( a \right) \sin \left( \pi\,a \right)  \left( a-b \right) }} \,,
\end{align}

which can be simplified using \eqref{WellKnown} to give

\begin{align} \label{C1A}
&\displaystyle {\mbox{$_3$F$_2$}(a,a,b;\,2\,a,b+1;\,1)}=-{\frac {\pi\,b\,\Gamma \left( 2\,a \right) {\mbox{$_3$F$_2$}(a,1-a,a-b;\,1,1+a-b;\,1)}}{ \Gamma \left( a \right) ^{2}\sin \left( \pi\,a \right)  \left( a-b \right) }}
\mbox{}\\ \nonumber
&+{\frac {{\pi}^{2}\,\Gamma \left( 2\,a \right) \Gamma \left( b+1 \right) }{  \Gamma \left( a \right)  ^{2}\sin \left( \pi\,a \right) \sin \left( \pi\, \left( a-b \right)  \right)  \Gamma \left( b+1-a \right)  ^{2}
\mbox{}\Gamma \left( 2\,a-b \right) }}\,.
\end{align}

Compare with \eqref{Sm1b} above and \eqref{hb2c} below. In contrast, it is also possible to operate on \eqref{Mix2} (after renaming the parameters) with the slightly different left-ordered operator sequence $ \frac{\partial}{\partial d},d=a,b=a, f\rightarrow b+1, e=2\,a$ to obtain a different representation


\begin{align} \label{hb2c}
&\displaystyle \frac {  \Gamma \left( a \right) ^{2}{_3F_2(a,a,b;\,2\,a,1+b;\,1)}}{\Gamma \left( 2\,a \right) b}
=-{\frac {\pi\,{_3F_2(a,1-a,a-b;\,1,1+a-b;\,1)}}{\sin \left( \pi\,a \right)  \left( a-b \right) }}
\mbox{} \\ \nonumber
&+{\frac {\sin \left( \pi \left( a-b \right)  \right) {_4F_3(a,1-a,a-b,a-b;\,1,1+a-b,1+a-b;\!1)}}{\sin \left( \pi \,b \right)  \left( a-b \right) ^{2}}}-
{\frac {\sin \left( \pi \left( a-b \right)  \right) \sin \left( \pi a \right) }{\pi\,\sin \left( \pi\,b \right) }{S_2(a,b)}}
\end{align}

where $S_{2}(a,b)$ is defined in \eqref{S2}. To complete the calculation, operate on each of the above with $\mathfrak{L}(b)$. In the case of \eqref{Hb3B}, the result simplifies to \eqref{Miller1Final}; in the case of \eqref{hb2c} the following representation is found:

\begin{align} \label{Hfb}
\displaystyle {\frac {h \left( a \right) }{{a}^{2}}}=&-{\frac {\Gamma \left( 2\,a \right)}{  \Gamma \left( a \right) ^{2}} S_{1}(a)}-2\, \pi\,\left( a-1 \right)\,{\frac { \Gamma \left( 2\,a+1 \right) {\mbox{$_5$F$_4$}(1,1,1,2-a,a+1;\,2,2,2,2;\,1)}}{\sin \left( \pi\,a \right)   \Gamma \left( a \right)  ^{2}}}
 \\ \nonumber
&+{\frac {2\,{\pi}^{2}\,\Gamma \left( 2\,a \right) }{  \Gamma \left( a \right) ^{2}
\sin \left( \pi\,a \right) } \left( \pi/6-{\frac { \left( \gamma+\Psi \left( a \right)  \right) \cos \left( \pi\,a \right) }{\sin \left( \pi\,a \right) }} \right) }
\end{align}

It is worth noting that 
\begin{align} \label{Note1}
&\displaystyle \sum _{k=1}^{\infty }{\frac {\Gamma \left( a+k \right) \Gamma \left( k+1-a \right) \Psi \left( k+1 \right) }{  k^{3}\,\Gamma \left( k \right)  ^{2}
}}-\sum _{k=1}^{\infty }{\frac {\Gamma \left( a+k \right) \Gamma \left( k+1-a \right) \Psi \left( k \right) }{k^{3} \,\Gamma \left( k \right) ^{2}}} \\ \nonumber
&=-{\frac {a \left(a -1 \right) \pi\,{\mbox{$_5$F$_4$}(1,1,1,1+a,2-a;\,2,2,2,2;\,1)}}{\sin \left( \pi\,a \right) }}
\end{align}

and that neither of the $_3F_2$ that appear in \eqref{FromB7} or \eqref{C1A} can be found in the database \cite{Milgram447}.

\subsection{Consider a two-part relation} 

Any $_3F_2(1)$ can be related to a maximum of nine other $_3F_2(1)$ by means of the two-part Thomae relations (see Appendix \ref{sec:AppendixC}). In the case of the particular $_3F_2(1)$ under consideration, we have four different possibilities

\begin{align} \label{Thomae1}
\displaystyle {\mbox{$_3$F$_2$}(a,a,b;\,2\,a,b+1;\,1)}&={\frac {{\mbox{$_3$F$_2$}(1,1,2\,a-b;\,1+a,1+a;\,1)}\Gamma \left( b+1 \right) \Gamma \left( 2\,a \right) }{\Gamma \left( b \right)  \left( \Gamma \left( 1+a \right)  \right) ^{2}
\mbox{}}} \\
&={\frac {{\mbox{$_3$F$_2$}(1,b-a+1,a;\,1+a,b+1;\,1)}\Gamma \left( 2\,a \right) }{\Gamma \left( a \right) \Gamma \left( 1+a \right) }} \\
&={\frac {{\mbox{$_3$F$_2$}(2\,a-b,a,a;\,1+a,2\,a;\,1)}\Gamma \left( b+1 \right) }{\Gamma \left( b-a+1 \right) \Gamma \left( 1+a \right) }} \\
&={\frac {{\mbox{$_3$F$_2$}(b-a+1,b-a+1,b;\,b+1,b+1;\,1)}\Gamma \left( 2\,a \right) }{\Gamma \left( 2\,a-b \right) \Gamma \left( b+1 \right), 
\mbox{}}}
\end{align}
none of which appear in the database \cite{Milgram447}, as expected; as well, each fails a test of the telescoping method \cite{Gosper}. However, each of these could be evaluated as a limiting case of one of the eight three-part transforms discussed in previous sections, giving 32 possibilities. I have performed this calculation using \eqref{Thomae1} and \eqref{Eq1p7b} as a basis, and found nothing useful beyond a long expression containing many intractable sums. A study of the remaining 31 possibilities is left as an exercise for the ambitious reader. 

\section{Reduce the problem}\label{sec:reduce}

When all else fails, it is sometimes helpful to consider a reduced problem - in this case consider the case $a\rightarrow n$. With reference to \eqref{Miller1Final}, the reduction to the Integers from the Reals will clearly require a limiting calculation because of the presence of the denominator factor $\sin(\pi a)$, so let $a=n+\epsilon$, and consider the limit $\epsilon \rightarrow 0$. The leading term of this limit on the right-hand side is of order $\epsilon^{-1}$, and when the (numerator) coefficient of this term is set to zero (for more detail, see the remark at the beginning of Section \ref{sec:d=n} below) we find

\begin{equation}
\mapleinline{inert}{2d}{hypergeom([1, 1, 1, 1+n, 2-n], [2, 2, 2, 2], 1) = 2*(gamma+Psi(n))^2/(n*(n-1))}{\[\displaystyle {\mbox{$_5$F$_4$}(1,1,1,1+n,2-n;\,2,2,2,2;\,1)}=2\,{\frac { \left( \gamma+\Psi \left( n \right)  \right) ^{2}}{n \left( n-1 \right) }}\]} \hspace{1cm} n>1\,.
\label{Ansm1a}
\end{equation}
For the general case of \eqref{Ansm1a}, see \eqref{5F4} below. For the case $n=1$ we have the well known definition
\begin{equation} \label{Z3}
_4F_3(1,1,1,1;2,2,2;1)=\zeta(3)\,.
\end{equation}

The next-to-leading term ($\epsilon^0$) produces the identity


\begin{align} \label{OldEq8}
&\displaystyle h(n)={\frac { \Gamma \left( 2\,n+1 \right)^{2}{\mbox{$_5$F$_4$}(1,1,n,n,2\,n;\,1+n,1+n,1+n,1+n;\,1)}}{  4\,\Gamma \left( 1+n \right)   ^{4}}}
+\frac { \left( -1 \right) ^{n}\,2\,n
\mbox{}\Gamma \left( 2\,n+1 \right) }{  \Gamma \left( n \right)  ^{2}} \\ \nonumber
&\hspace{1cm}\times  \bigg{(}  \left( \gamma+\Psi \left( n \right)  \right) \Psi^{\prime} \left(n \right)
+\frac{1}{4}\,\sum _{k=0}^{n-2}{\frac { \left( -1 \right) ^{k}\Gamma \left( n+k+1 \right)  \left( \Psi \left( n-1-k \right) -\Psi \left( n+k+1 \right)  \right) }{  \Gamma \left( k+1 \right)   ^{2}
\mbox{}\Gamma \left( n-1-k \right)  \left( k+1 \right) ^{4}}} \bigg{)}\,.
\end{align}

Alternatively, from \cite[Entry 29]{Milgram447} and using the approach based upon the application of $\mathfrak{L}(b)$ to a $_3F_2(1)$, after taking the appropriate limit, we obtain


\begin{align}
&\displaystyle {\mbox{$_3$F$_2$}(n,n,b;\,2\,n,b+1;\,1)}=\frac {\Gamma \left( 2\,n \right) }{\Gamma \left( 2\,n-b \right) } \left( {\frac { \left( \pi\,\cot \left( \pi\,b \right) -2\,\Psi \left( n \right) +2\,\Psi \left( -n+b+1 \right)  \right)   \Gamma \left( n-b \right) ^{2} }{\Gamma \left( -b \right)  \Gamma \left( n \right)  ^{2}}}
\right. \\ \nonumber
&\left.  +b \left( -1 \right) ^{n} \left( \sum _{k=0}^{n-1}{\frac {\Gamma \left( n+k-b \right)  \left( -1 \right) ^{k}\Psi \left( n-k \right) }{\Gamma \left( k+1 \right)  ^{2}\Gamma \left( n-k \right) }}-\sum _{k=0}^{n-1}{\frac {\Gamma \left( n+k-b \right)  \left( -1 \right) ^{k}\Psi \left( n+k-b \right) }{\Gamma \left( k+1 \right) ^{2}\Gamma \left( n-k \right) }}
\mbox{} \right)  \right)\,. 
\end{align}

The second sum can be evaluated from \eqref{CxSum} to eventually yield

\begin{align} \label{Cx0b}
&\displaystyle {\mbox{$_3$F$_2$}(n,n,b;\,2\,n,b+1;\,1)}=\frac {\Gamma \left( 2\,n \right) }{\Gamma \left( 2\,n-b \right) } \\ \nonumber
&\times\left( b \left( -1 \right) ^{n}\sum _{k=0}^{n-1}{\frac {\Gamma \left( n+k-b \right)  \left( -1 \right) ^{k}\Psi \left( n-k \right) }{  \Gamma \left( k+1 \right)^{2}\Gamma \left( n-k \right) 
}}+{\frac { \left( -2\,\Psi \left( n \right) +\Psi \left( b \right)  \right)   \Gamma \left( n-b \right)  ^{2}}{\Gamma \left( -b \right)   \Gamma \left( n \right)   ^{2}}}
\mbox{} \right)\,. 
\end{align}
After operating on \eqref{Cx0b} with $\mathfrak{L}(b)$ (see \eqref{MathL}) and setting $b=n$, we ultimately obtain (also see \eqref{Eq82}) a second representation


\begin{align} \label{Cx0d}
&\displaystyle \frac{h \left( n \right)}{{n}^{2}} = \Gamma \left( 2\,n \right)  \left( -1 \right) ^{n} \left( {\frac {1}{\Gamma \left( n \right) 
\mbox{}}\sum _{k=1}^{n-1}{\frac { \left( -1 \right) ^{k}\Psi \left( n-k \right)  \left( \Psi \left( k \right) -\Psi \left( n \right)  \right) }{{k}^{2}\Gamma \left( n-k \right) \Gamma \left( k \right) }}} \right. \\ \nonumber
&\left.+{\frac {\Psi \left( n \right) {\pi}^{2}/12-3/2\,\Psi \left( n \right) ^{3}-3\,  \Psi \left( n \right)^{2}\gamma+3/2\, \left( -{\gamma}^{2}+\Psi^{\prime} \left(n \right)  \right) \Psi \left( n \right) 
\mbox{}+2\,\Psi^{\prime} \left(n \right) \gamma+\Psi^{(2)} \left( n \right)/2 }{ \Gamma \left( n \right)^{2}}} \right)\,.
\end{align}

For the sake of completeness, in the case that $n\rightarrow -n+\epsilon$, from \eqref{Miller1Final}  we find the first three terms of the expansion in $\epsilon$:

\begin{align} \label{Hexp}
&\displaystyle h(-n+\epsilon) \approx {\frac {n \left( -1 \right) ^{n}  \Gamma \left( n+1 \right)  ^{2}}{2\,\Gamma \left( 2\,n \right) {\epsilon}^{2}}}-{\frac { \Gamma \left( n+1 \right)  ^{2} \left( -1 \right) ^{n} \left( 2\,n\gamma-2\,n\Psi \left( 2\,n \right) +4\,n\Psi \left( n+1 \right) +1 \right)
\mbox{}}{2\,\Gamma \left( 2\,n \right) \epsilon}} \\ \nonumber
&+\frac{\left( -1 \right) ^{n} \Gamma \left( n+1 \right)   ^{2}}{\Gamma \left( 2\,n \right)} \left(\gamma-\Psi \left( 2\,n \right)+2\,\Psi \left( n+1 \right) 
  + n \bigg{(} 2\,\gamma \left(\Psi \left( n+1 \right) -\Psi \left( 2\,n \right)   \right) + \Psi \left( 2\,n \right) ^{2} \right.  \\ \nonumber
& \left. -4\,\Psi \left( n+1 \right) \Psi \left( 2\,n \right) +3\, \Psi \left( n+1 \right)^{2}
\mbox{}-{\pi}^{2}/12\,+3/2\,\Psi^{\prime} \left(n+1 \right) -\Psi^{\prime} \left(2\,n \right)  \bigg{)} \right)+\dots  
\end{align}

\section{The method of ``Hail Mary"}\label{sec:hail}

Since none of the various methods attempted have thus far yielded a closed form for \eqref{question}, there is one ``last-ditch" method available, which relies more on experience and luck than mathematical rigour. It is included here for both its educational value and because, in this case, it does produce a different representation from what has been obtained so far. Analysts who do not have to solve practical problems might wish to stop reading this section at this point. \newline

Near the end of Section \ref{sec:d=n} will be found a discussion outlining how it is sometimes beneficial to change an expression that involves only integer variables into equivalent expressions involving real variables by the simple expedient of changing $n=a$ and checking the result numerically, provided that the modified expression makes sense after the change. Sometimes this procedure works, and possible reasons are given; many times it does not - the presence of terms involving $\sin(\pi a)$ in a comparison of \eqref{Ansm1a} and \eqref{5F4} exemplifies why. In this section a similar change will be performed, despite the fact that the modified expression does not make sense after the change.\newline

Consider the representation \eqref{Cx0d} and replace $n=a$. The term $(-1)^n$ is reasonably transformed into $\exp(i\pi a)$, and all other terms are well-defined functions of a complex variable, but what are we to make of the summation with a non-integer limit? In both the derivation of \eqref{Eq1p7b} and in a later section (see the latter part of Section \ref{sec:d=n}) this action was employed and justified - simply extend any sums to infinity, relying on the presence of denominator gamma functions ($\Gamma(n-k)$ in this case) to effectively terminate the sum. That action does not bode well here, because of the presence of the digamma function in the numerator of \eqref{Cx0d}; that is $\Psi(n-k)/\Gamma(n-k)$ is finite for $k>n$. Thus new terms are added if the sum were to be extended to infinity by fiat. In this regard, M\"uller and Schleicher have studied that question and provide a well-reasoned prescription \cite[Eq. (4)]{AddTerms} for any indexed function $f(k)$, imposing simple conditions on the asymptotics, as follows:  

\begin{equation}
\sum_{k=x}^{y}f(k)=\sum_{k=1}^{\infty}(f(k+x-1)-f(k+y))\,.
\label{NonInt}
\end{equation}
After replacing $n:=a$ and applying \eqref{NonInt} to the sum in \eqref{Cx0d} we obtain a fairly long and complicated complex expression, whose imaginary part must vanish if we temporarily restrict $a \in \mathbb{R}$ (e.g. see Appendix \ref{sec:MeijerG}). This procedure yields the following identity



\begin{align} \label{AnsIm}
&\displaystyle \sum _{k=1}^{\infty }{\frac { \left( -1 \right) ^{k}\Psi \left( a-k \right)  \left( \Psi \left( k \right) -\Psi \left( a \right)  \right) }{{k}^{2}\Gamma \left( k \right) \Gamma \left( a-k \right) 
\mbox{}}} =\frac{1}{\Gamma \left( a \right)} 
\bigg{(} 3/2\,  \Psi \left( a \right) ^{3}+3\,\gamma\,  \Psi \left( a \right) ^{2}  \\ \nonumber
&\hspace{1cm}+ \left( -3/2\,\Psi^{\prime} \left(a \right) +3/2\,{\gamma}^{2}-1/12\,{\pi}^{2} \right) \Psi \left( a \right) 
\mbox{}-2\,\gamma\,\Psi^{\prime} \left(a \right) -1/2\,\Psi^{(2)} \left(a \right) { }\bigg{)}
\end{align}

which, mirabile dictu, is numerically satisfied for $\Re(a)>1$, where the sum converges. To further check \eqref{AnsIm}, see \eqref{AnsImG} where it is explored in the reverse limit $a\rightarrow n$, yielding a representation for what may be a new generalized Euler sum, verifiable in some special cases. Furthermore, with reference to \eqref{Lemma6}, \eqref{AnsIm} can also be rewritten in the form of a possibly new sum 


\begin{align} \label{AnsImb}
\displaystyle \sum _{k=1}^{\infty }{\frac { \left( -1 \right) ^{k}\Psi \left( a-k \right) \Psi \left( k \right) }{{k}^{2}\Gamma \left( k \right) \Gamma \left( a-k \right) }}&={\frac { \Psi \left( a \right)  ^{3}}{2\,\Gamma \left( a \right) }}
\mbox{}+{\frac {2\,\gamma\, \Psi \left( a \right) ^{2}}{\Gamma \left( a \right) }}-{\frac { \left( {\pi}^{2}-18\,{\gamma}^{2}+6\,\Psi^{\prime} \left(a \right)  \right) \Psi \left( a \right) }{12\,\Gamma \left( a \right) }}\\ \nonumber
&-{\frac {4\,\gamma\,\Psi^{\prime}\left(a \right) +\Psi^{(2)} \left(a \right) }{2\,\Gamma \left( a \right) }}\,.
\end{align}

As noted elsewhere, the right-hand side of \eqref{AnsImb} represents the left-hand side for arbitrary values of $a$ by the principle of analytic continuation. However, even though the sum converges for $a>1$, it will not truncate to a finite sum when $a=n$ because of the numerator factor $\Psi(n-k)$, and thus cannot be used as a bootstrap to further evaluate \eqref{Cx0d}.\newline 

Finally, with the help of \eqref{Lemma2} and \eqref{Lemma5} followed by the conversion of the remaining sum into hypergeometric form, the real part of \eqref{Cx0d} becomes

\begin{align} \label{AnsImB}
\displaystyle {\frac {h \left( a \right) }{{a}^{2}}}&=-{\frac {\Gamma \left( a+1/2 \right) {4}^{a}{\mbox{$_4$F$_3$}(1,1,a,a;\,1+a,1+a,1+a;\,1)}}{ \sqrt{\pi}\Gamma \left( 1+a \right) {a}^{2}}}
\mbox{}\\ \nonumber
&-{\frac { \left( \Psi^{(2)} \left(a/2 \right) -\Psi^{(2)} \left(a/2+1/2 \right)  \right) \Gamma \left( a+1/2 \right)
\mbox{}{4}^{a}}{8\, \sqrt{\pi}\Gamma \left( a \right) }}
\end{align} 

a representation that outwardly differs from the other representations obtained up to this point. However, once established, however shakily, \eqref{AnsImB} can be shown analytically to be equivalent to the other representations; it is numerically verifiable and therefore appears to be true, despite its unusual lineage. In the limiting case $a=1$, \eqref{AnsImB} consistently reduces to \eqref{Z3}, i.e. $h(1)=\zeta(3)$. \newline

\section{(Unforeseen) Consequences} \label{sec:Consq}

\subsection{Proof of some new $_4F_3(1)$ and $_5F_6(1)$}

An unexpected result materializes when \eqref{Miller2Final} is compared to \eqref{Miller2d}. From that comparison, a closed form of the following (non-terminating, Saalsch\"utzian) $_4F_3(1)$ is found. It appears to be new:
\begin{equation}
\mapleinline{inert}{2d}{hypergeom([1, 1, 2-a, a+1], [2, 2, 2], 1) = (2*(Psi(a)+gamma))/(a*(-1+a))+(1/2)*sin(Pi*a)*(Psi^{\prime}( (1/2)*a)-Psi^{\prime}( (1/2)*a+1/2))/(a*(-1+a)*Pi)}{\[\displaystyle {\mbox{$_4$F$_3$}(1,1,2-a,a+1;\,2,2,2;\,1)}=2\,{\frac {\Psi \left( a \right) +\gamma}{a \left(a -1 \right) }}+{\frac {\sin \left( \pi\,a \right)  \left( \Psi^{\prime} \left(a/2 \right) -\Psi^{\prime} \left( a/2+1/2 \right)  \right) }{2\,a\, \pi\, \left(a -1 \right) }}\]} \,.
\label{4F3}
\end{equation}
Inspired by this observation, if the following series of operations, ordered from the left,
\begin{equation} \label{reduction}
e=2,f=2,g=2,a\rightarrow b,b\rightarrow 1
\end{equation}

are applied to \eqref{Miller1} the following transformation (an apparent generalization of \eqref{4F3}) results

\begin{align} \label{Y1b1}
&\displaystyle {\mbox{$_4$F$_3$}(1,1,c,d;\,2,2,2;\,1)}=-{\frac {\sin \left( \pi\,c \right) \Gamma \left( c \right) {\mbox{$_4$F$_3$}(1,c,c,c;\,2,c+1,2+c-d;\,1)}\Gamma \left( 1-d \right) }{c\Gamma \left( 2+c-d \right) 
\mbox{}\pi}} \\ \nonumber
&-{\frac {\Psi \left( -1+c \right) }{ \left( -1+c \right)  \left( -1+d \right) }}-{\frac {\Psi \left( 2-d \right) }{ \left( -1+c \right)  \left( -1+d \right) }}
\mbox{}-{\frac {\sin \left( \pi\,c \right) \Gamma \left( c \right) \Gamma \left( 1-d \right) }{\Gamma \left( 1+c-d \right)  \left( -1+c \right) ^{3}\pi}}-2\,{\frac {\gamma}{ \left( -1+c \right)  \left( -1+d \right) }}\,.
\end{align}

The $_4F_3(1)$ appearing on the right-hand side of this equation can be reduced using \eqref{WellKnown}, and when this is done, \eqref{Y1b1} becomes


\begin{align} \label{Y1b2}
\displaystyle {\mbox{$_4$F$_3$}(1,1,c,d;\,2,2,2;\,1)}=&-{\frac {\sin \left( \pi\,c \right) \Gamma \left( c \right) \Gamma \left( 1-d \right) {\mbox{$_3$F$_2$}(c-1,c-1,c-1;\,c,1+c-d;\,1)}}{\Gamma \left( 1+c-d \right) 
\mbox{} \left( -1+c \right) ^{3}\pi}}\\ \nonumber
&-{\frac {2\,\gamma+\Psi \left( c-1 \right) +\Psi \left( 2-d \right) }{ \left( c-1 \right)  \left( d-1 \right) }}\,.
\end{align}

For many choices of the parameter $d$, the $_3F_2(1)$ on the right-hand side is known as a special case of, or contiguous to, Watson's theorem \cite[Entry 1]{Milgram447}; in the following sub-sections, this will be explored. The case $d\rightarrow 1$ reduces the left-hand side of \eqref{Y1b1} to the known result \eqref{Prud9} below. The case $c=d$ yields a result equivalent to $a\rightarrow c =1$ in \eqref{MStep2}. {\bf Remark:} attempting to apply \eqref{WellKnown} to the left-hand side of \eqref{Y1b1} or \eqref{Y1b2} simply leads to a trivial identity.
\subsubsection{Special case d=c+1-a}

Because the left-hand side of \eqref{Y1b2} is symmetrical under the exchange $c\leftrightarrow d$ and the right-hand side is not, equating the two right-hand sides after the interchange gives rise to the transformation


\begin{align} \label{H12}
\displaystyle {\mbox{$_3$F$_2$}} \left( c,c,c;c+1,1+a;1 \right)& ={\frac {{\pi}^{2}c\Gamma \left( 1+a \right) {\mbox{$_3$F$_2$}} \left( c-a,c-a,c-a;c+1-a,1-a;1 \right) }{ \left( c-a \right) ^{3} \Gamma \left( -c+a \right)  ^{2}
\mbox{} \sin^2 \left( \pi\,c \right) \Gamma \left( c \right)  ^{2}\Gamma \left( 1-a \right) }}\\ \nonumber
&-{\frac {\pi\,\sin \left( \pi\,a \right) \Gamma \left( c-a \right) c\Gamma \left( 1+a \right) }{ \sin^2 \left( \pi\,c \right) \Gamma \left( c \right) }}
\end{align}
being almost certainly a special case of the three-part transformations discussed in Appendix \ref{sec:AppendixC}. Of interest are the two-part (Thomae) equivalents of the $_3F_2$ that appears on the right hand-side of \eqref{H12} giving

\begin{align} \nonumber
\displaystyle \mbox{$_3$F$_2$}&(c,c,c;\,c+1,1+a;\,1)\\ \nonumber
&={\frac {c\,\Gamma \left( 2-2\,c+a \right) \Gamma \left( 1+a \right) {\mbox{$_3$F$_2$}(1,1-c,2-2\,c+a;\,2-c,2-c;\,1)}}{ \left( 1-c \right) ^{2} \Gamma \left( 1-c+a \right)  ^{2}}} \\ 
&\hspace{3cm}-{\frac {c\,\sin \left( \pi\,a \right) \Gamma \left( 1-c \right) \Gamma \left( c-a \right) \Gamma \left( 1+a \right) }{\sin \left( \pi\,c \right) }}
\label{t1a}\\  \nonumber
&=-{\frac {\pi\,c\Gamma \left( c-a \right) \Gamma \left( 1-c \right) a\Gamma \left( 2-2\,c+a \right) {\mbox{$_3$F$_2$}(1-c,1-c,c-a;\,2-c,1-a;\,1)}}{ \left( c-1 \right) 
\mbox{}\sin \left( \pi\,a \right)  \left( \Gamma \left( 1-a \right)  \right) ^{2} \left( \Gamma \left( 1-c+a \right)  \right) ^{2}}}\\
&\hspace{3cm}-{\frac {\pi\,c\Gamma \left( c-a \right) \Gamma \left( 1-c \right) a}{\Gamma \left( 1-a \right) \sin \left( \pi\,c \right) }}
\label{t3} \\ 
&=\displaystyle {\frac {c\Gamma \left( 1+a \right) \Gamma \left( 2-2\,c+a \right) {\mbox{$_3$F$_2$}(1,1,c-a;\,2-c,c+1-a;\,1)}}{ \left( c-a \right)  \Gamma \left( 1-c+a \right)^{2}
\mbox{} \left( 1-c \right) }}-{\frac {\pi\,c\,a\,\Gamma \left( c-a \right) \Gamma \left( 1-c \right) }{\Gamma \left( 1-a \right) \sin \left( \pi\,c \right) }}
\label{t4} \\ \nonumber
&=\displaystyle {\frac {c\,a \sin^2 \left( \pi\, \left( c-a \right)  \right)\Gamma \left( 1-c \right)\Gamma \left( c-a \right)  ^{3}
\mbox{}\Gamma \left( 2-2\,c+a \right) {\mbox{$_3$F$_2$}(1-c,1-c,c-a;\,2-c,1-a;\,1)}}{\pi\, \left( 1-c \right) \sin \left( \pi\,a \right) \Gamma \left( 1-a \right) ^{2}}}\\
&\hspace{3cm}-{\frac {\pi\,c\,a\,\Gamma \left( c-a \right) \Gamma \left( 1-c \right) }{\Gamma \left( 1-a \right) \sin \left( \pi\,c \right) }}\,.
\label{t6}
\end{align}

None of the above are contained in the database \cite{Milgram447} in their general form. Note that the validity of \eqref{t1a} is limited to $2c-1<a<c+2$, along with $c<3$ and that of \eqref{t6} requires $c<1$ and $a>2c-2$, thereby defining an overlap region that establishes the analytic continuation of either side with reference to the other. However, in both cases \eqref{t3} and \eqref{t6} the right-hand sides have parametric excess equal to unity, in which case, the terminating instance of either is Saalsch\"utzian and therefore summable. Thus it is tempting to investigate both of these transformations corresponding to the special cases $c=1+n$ or $a=c+n$. However, this fails, because the Saalsch\"utz formula does not apply in a limiting sense, and in both cases the application of Saalsch\"utz' formula leads to representations in which limits must be further employed. So, rather than using a valid, but inappropriate (Saalsch\"utzian) equality, consider:
\begin{itemize}
\item

\eqref{t3} in the limit $a\rightarrow c+n$. After considerable calculation, the result is

\begin{align} \label{T0}
&\displaystyle {\frac {{\mbox{$_3$F$_2$}(c,c,c;\,c+1,1+c+n;\,1)}}{\Gamma \left( 1+c+n \right) }}={\frac { \left( \Psi \left( n+1 \right) -\Psi \left( 2-c+n \right)  \right) c\pi}{\Gamma \left( n+1 \right) \Gamma \left( c \right) \sin \left( \pi\,c \right) }}\\ \nonumber 
&-{\frac {c \Gamma \left( 2-c+n \right)   ^{2}{\mbox{$_4$F$_3$}(1,1,2-c+n,2-c+n;\,2-c,n+2,3-c+n;\,1)}}{ \left( -2+c-n \right) \Gamma \left( n+2 \right) 
\mbox{}\Gamma \left( 2-c \right)  \Gamma \left( n+1 \right) ^{2}}} \\ \nonumber
&-{\frac {c \left( -1 \right) ^{n}\Gamma \left( 2-c+n \right) }{ \left( \Gamma \left( n+1 \right)  \right) ^{2}}\sum _{k=0}^{n}{\frac { \left( -1 \right) ^{k}\Gamma \left( 1-c+k \right)  \left( \Psi \left( n-k+1 \right) -\Psi \left( c+n-k \right)  \right) }{ \left( c-k-1 \right) \Gamma \left( 1-c-n+k \right) 
\mbox{}\Gamma \left( n-k+1 \right) \Gamma \left( k+1 \right) }}}\,.
\end{align} 
In comparison, the specific form of the left-hand side of \eqref{t3} when $a=c+n$ is easily established (\cite[Entry 27]{Milgram447}, that is 

\begin{align}\label{BC0}
\displaystyle {\mbox{$_3$F$_2$}(c,c,c;\,c+1,1+c+n;\,1)}&=-{\frac {c\Gamma \left( 1-c \right) \Psi \left( c \right) \Gamma \left( 1+c+n \right) }{\Gamma \left( n+1 \right) }}\\ \nonumber
&-{\frac {c\,\Gamma \left( 2-c+n \right) \Gamma \left( 1+c+n \right) }{\Gamma \left( n+1 \right) }\sum _{k=0}^{n}{\frac { \left( -1 \right) ^{k}\Psi \left( k+1 \right) }{ \left( c-k-1 \right) \Gamma \left( n-k+1 \right) \Gamma \left( k+1 \right) 
\mbox{}}}}
\end{align}
and so, by comparing the right-hand sides of \eqref{T0} and \eqref{BC0} and redefining the variable $c$, we find a possibly new transformation

\begin{align} \label{BC1b}
&\displaystyle {\frac {{\mbox{$_4$F$_3$}(1,1,c,c;\,c-n,2+n,c+1;\,1)}}{ c\,\Gamma \left( 2+n \right) \Gamma \left( c-n \right)\Gamma \left( 1+n \right) }}={\frac { \left( -\Psi \left( -c+n+2 \right) 
 -\Psi \left( 1+n \right) +\Psi \left( c \right)  \right) \Gamma \left( -1+c-n \right) }{ \Gamma \left( c \right) ^{2}}} \\ \nonumber
& \hspace{2cm}-{\frac {1}{\Gamma \left( c \right) }\sum _{k=0}^{n}{\frac { \left( -1 \right) ^{k}\Psi \left( k+1 \right) }{ \left( n-k+1-c \right) \Gamma \left( n-k+1 \right) \Gamma \left( k+1 \right) }}}\\ \nonumber 
&\hspace{2cm}-{\frac { \left( -1 \right) ^{n}}{\Gamma \left( c \right) \Gamma \left( 1+n \right) }\sum _{k=0}^{n}{\frac { \left( \Psi \left( 2-c+k+n \right) -\Psi \left( k+1 \right)  \right)  \left( -1 \right) ^{k}\Gamma \left( 2-c+k+n \right) }{\Gamma \left( 2-c+k \right) 
\mbox{} \left( 1-c+k \right) \Gamma \left( k+1 \right) \Gamma \left( n-k+1 \right) }}}
\end{align}\newline

\item
\eqref{t3} with $c=n+1+\epsilon$ in the limit $\epsilon\rightarrow 0$. As discussed elsewhere (see remark at the start of Section \ref{sec:d=n}), a leading (divergent) term of order $\epsilon^{-1}$ arises, and when the numerator coefficient of that term is set to zero, we find (after redefining $a:=a+n-1,n:=n+1$) the closed form identity

\begin{equation}
\mapleinline{inert}{2d}{T3Cbm1 := hypergeom([1, 1, a, -n], [2, 2, -n+a-1], 1) = (gamma+Psi(2+n)+Psi(a-1)-Psi(a-2-n))*(2+n-a)/((1+n)*(1-a))}{\[\displaystyle {\mbox{$_4$F$_3$}(1,1,a,-n;\,2,2,a-n-\!1;1)}={\frac { \left( \gamma+\Psi \left( 2+n \right) +\Psi \left( a-1 \right) -\Psi \left( a-2-n \right)  \right) 
\mbox{} \left( 2+n-a \right) }{ \left( 1+n \right)  \left(1-a \right) }}\]}.
\label{T3Cbm1}
\end{equation}

Since the left-hand side of \eqref{T3Cbm1} is terminating and Saalsch\"ultzian, it is summable and in that sense \eqref{T3Cbm1} is possibly not new, but represents a limiting case that could possibly have been obtained by other means.\newline

The left-hand side of \eqref{t3} can also be reduced using \cite[Entry 26]{Milgram447} to eventually obtain
\end{itemize}


\begin{align} \label{Ht0}
&\displaystyle {\frac {\Gamma \left( 1+n \right) \Gamma \left( a \right) {\mbox{$_3$F$_2$}(1+n,1+n,1+n;\,2+n,1+a+n;\,1)}}{\Gamma \left( 1+a+n \right)  \left( 1+n \right) }} \\ \nonumber
&=\Gamma \left( 1+n \right) \Gamma \left(a -n \right) \sum _{k=0}^{n-1}{\frac { \left( -1 \right) ^{k} \left( \Psi \left( a-k \right) +\Psi \left( n-k \right) -\Psi \left( k+1 \right)  \right) }{ \left( n-k \right) \Gamma \left( a-k \right) \Gamma \left( k+1 \right) }} \\ \nonumber
&-\frac{(-1)^n}{2}\, \left(  \left( \Psi \left( a-n \right) -\Psi \left( a \right)  \right)  \left( 3\,\Psi \left( a-n \right) -2\,\gamma-2\,\Psi \left( 1+n \right) -\Psi \left( a \right)  \right) -3\,\Psi^{\prime} \left(a-n \right) 
\mbox{}+\Psi^{\prime} \left(a \right)  \right).  
\end{align}
Now compare the limiting case \eqref{BC0} using $c\rightarrow n+1$ and \eqref{Ht0} with $a=n+1$, and after some simplification, obtain what may possibly be a new sum:

\begin{equation}
\mapleinline{inert}{2d}{HTxA := Sum((-1)^k*Psi(k+1)/((k+1)^2*GAMMA(n-k)*GAMMA(k+1)), k = 0 .. -1+n) = (1/12)*(Pi^2-18*gamma^2-24*gamma*Psi(1+n)-6*Psi(1+n)^2-6*Psi(1, 1+n))/GAMMA(1+n)}{\[\displaystyle\sum _{k=0}^{n-1}{\frac { \left( -1 \right) ^{k}\Psi \left( k+1 \right) }{ \left( k+1 \right) ^{2}\Gamma \left( n-k \right) \Gamma \left( k+1 \right) 
\mbox{}}}={\frac {{\pi}^{2}-18\,{\gamma}^{2}-24\,\gamma\,\Psi \left( 1+n \right) -6\, \Psi \left( 1+n \right)^{2}-6\,\Psi^{\prime} \left(1+n \right) }{12\,\Gamma \left( 1+n \right) }}\]},
\label{HTxA}
\end{equation}

after which \eqref{BC0} simplifies to

\begin{align} \label{HTx0A}
\displaystyle & \frac {\Gamma \left( 1+n \right) ^{2}{\mbox{$_3$F$_2$}(1+n,1+n,1+n;\,2+n,2+2\,n;\,1)}}{\Gamma \left( 2+2\,n \right)  \left( 1+n \right) }
=-\Gamma \left( 1+n \right) \sum _{k=0}^{n-1}{\frac { \left( -1 \right) ^{k}\Psi \left( k+1 \right) }{ \left( n-k \right) ^{2}\Gamma \left( n-k \right) \Gamma \left( k+1 \right) }} \\ \nonumber
&\hspace{5cm}- \left( \gamma\,\Psi \left( 1+n \right) + \Psi \left( 1+n \right) ^{2}-\Psi^{\prime} \left( 1+n \right)  \right)  \left( -1 \right) ^{n}\,.
\end{align}
Notice the symmetry of the components of the sums in \eqref{Ht0} under reversal when $a=n+1$.

\subsubsection{Special case $d=n$} \label{sec:d=n}

\hspace{2cm}{\bf Case: $n>1$}\newline

({\bf Remark:} The following focusses on another occurrence in this paper of a serendipitous condition that, when it arises and when it is recognized, sometimes leads to the discovery of new identities. Therefore, it is presented in more detail than would otherwise be warranted.) Consider the limit $d\rightarrow n$ applied to \eqref{Y1b2}. This is easiest done by setting $d=n+\epsilon$ and evaluating the limit $\epsilon \rightarrow 0$ to the series representation of the right-hand side. When this calculation is performed, it is found that the leading term is of order $1/\epsilon$. Since the left-hand side is obviously finite with respect to $\epsilon$ in this limit, provided that $n<\left \lfloor 4-c \right \rfloor$, the only possibility is that the numerator of this divergent term must be identically zero. Indeed, after the numerator expression is carefully calculated, equated to zero, simplified and reordered, the following, tentatively new, identity arises
\begin{align} \label{Clm1A}
\displaystyle {\mbox{$_4$F$_3$}(1,c,c,c;\,2,c+1,2+c-n;\,1)}&=-\frac {c\,\Gamma \left( 2+c-n \right) }{c-1} \\ \nonumber
&\times\left(  \left( -1 \right) ^{n}\Gamma \left( 1-c \right) \Gamma \left( n-1 \right) +{\frac {1}{ \left( c-1 \right) ^{2}\Gamma \left( c-n+1 \right) }}
\mbox{} \right) \,.
\end{align} 
When \eqref{WellKnown} is applied, \eqref{Clm1A} reduces to the following identity 
\begin{equation}
\mapleinline{inert}{2d}{hypergeom([c, c, c], [c+1, 2+c-n], 1) = c*(-1)^n*GAMMA(1-c)*GAMMA(n-1)*GAMMA(2+c-n)}{\[\displaystyle {\mbox{$_3$F$_2$}(c,c,c;\,c+1,2+c-n;\,1)}=c \left( -1 \right) ^{n}\Gamma \left( 1-c \right) \Gamma \left( n-1 \right) \Gamma \left( 2+c-n \right) \]}\,,
\label{Clm1B}
\end{equation}
which corresponds to a special {limiting case} of entry 35 of \cite{Milgram447}, corresponding to a limiting case of \cite[Eq. (7.4.4.15)]{prudnikov}. Thus, interpreted as a limiting case, \eqref{Clm1B} is not new. Alternatively, since a top parameter exceeds a bottom one by a positive integer if $n>2$, either of the above are computable from the Minton-Karlsson algorithm \cite{Minton},\cite{Karlsson} and Appendix \ref{sec:Misc}.  Working backwards from \eqref{Clm1B} suggests that \eqref{Clm1A} is also ``not new" in the same sense, but to the best of my knowledge, it does not appear explicitly in the literature. New or not, this sequence of calculations demonstrates how unexpected identities can sometimes arise.\newline

Reverting to a discussion of the leading term of order $\epsilon^{0}$ in the case $d=n+\epsilon$, we find the (possibly new) object of the exercise

\begin{align} \label{Clm0}
\displaystyle {\mbox{$_4$F$_3$}(1,1,c,n;\,2,2,2;\,1)}&=-{\frac {\sin \left( \pi\,c \right)  \left( -1 \right) ^{n}}{\Gamma \left( c \right) \pi\,\Gamma \left( n \right) }\sum _{k=0}^{\infty }{\frac { \Gamma \left( c+k-1 \right)^{2}\Psi \left( 1+c-n+k \right) }{ \left( c+k-1 \right) \Gamma \left( k+1 \right) 
\mbox{}\Gamma \left( 1+c-n+k \right) }}} \\ \nonumber
& -{\frac {2\,\Psi \left( n \right) +\Psi \left( c \right) +2\,\gamma}{ \left( c-1 \right)  \left( n-1 \right) }}
\mbox{}+{\frac {2\,c+n-3}{ \left( n-1 \right) ^{2} \left( c-1 \right) ^{2}}}\,.
\end{align}

Both sides of \eqref{Clm0} are convergent if $c+n<4$.\newline

\hspace{2cm}{\bf Case: $n=1$}\newline

In the case that $n=1$, it is clear that a limiting ($d\rightarrow1$) operation will be required, and when this is performed as in the previous subsection, a leading term of order $\epsilon^{-1}$ arises. Setting the numerator to zero yields 
\begin{equation}
\mapleinline{inert}{2d}{hypergeom([c, c, c], [c+1, c+1], 1) = Pi*c^2*(Psi(c)+gamma)/sin(Pi*(c+1))}{\[\displaystyle {\mbox{$_3$F$_2$}(c,c,c;\,c+1,c+1;\,1)}=-{\frac {\pi\,{c}^{2} \left( \Psi \left( c \right) +\gamma \right) }{\sin \left( \pi\, c \right) }}\]}\,,
\label{Hccc}
\end{equation}

a known result corresponding to \cite[Entry 27]{Milgram447}. Considering the term of order $\epsilon^{0}$, and applying \eqref{Hccc} leads to


\begin{align} \label{D1}
\displaystyle {\mbox{$_4$F$_3$}(1,1,1,c;\,2,2,2;\,1)}=&{\frac {\sin \left( \pi\,c \right) }{ \pi\, \Gamma \left( c \right) }\sum _{k=0}^{\infty }{\frac { \Gamma \left( c+k-1 \right)\Psi \left( c+k \right) }{\left( c+k-1 \right) ^{2}\Gamma \left( k+1 \right) 
}}} \\ \nonumber
&+{\frac {\gamma\, \left( \Psi \left( c-1 \right) +\gamma \right) }{c-1}}
\mbox{}+{\frac {{\pi}^{2}}{6\,(c-1)}}\,.
\end{align}

Helpfully, Prudnikov et. al. \cite[Eq. 7.5.3(9)]{prudnikov}, (also \eqref{pFq}) have, by other means, listed 

\begin{align}\label{Prud9}
\displaystyle \mbox{$_4$F$_3$}(a,a,a,b;&\,1+a,1+a,1+a;\,1)=\\ \nonumber
&\frac {{a}^{3}\Gamma \left( a \right) \Gamma \left( 1-b \right)  \left( \Psi^{\prime} \left(a \right) -\Psi^{\prime} \left( 1+a-b \right) + \left( \Psi \left( a \right) -\Psi \left( 1+a-b \right)  \right) ^{2} \right) 
\mbox{}}{2\,\Gamma \left( 1+a-b \right) }
\end{align}

so, by applying a simple change $(a=1, b=c)$ to \eqref{Prud9} and comparing the right-hand sides of  \eqref{D1} and modified  \eqref{Prud9} we find what appears to be a new identity

\begin{align} \label{D1S}
&\displaystyle \sum _{k=0}^{\infty }{\frac { \Gamma \left( c+k-1 \right) \Psi \left( c+k \right) }{ \left( c+k-1 \right) ^{2}\Gamma \left( k+1 \right) }}=-{\frac {\Gamma \left( c-1 \right) {\pi}^{3}}{4\,\sin \left( \pi\,c \right) }}\\ \nonumber
&\hspace{1cm}-\frac{\pi\,\Gamma \left( c-1 \right)}{\sin \left( \pi\,c \right)}\, { \left( \frac{1}{2} \left( \gamma+\Psi \left( 2-c \right)  \right) ^{2}+{\gamma}^{2}+\gamma\,\Psi \left( c-1 \right) -\frac{1}{2}\Psi^{\prime} \left(2-c \right)  \right)  }\,.
\end{align}

Further interesting results can be obtained from the above ($d=1$) case by considering the substitution $c= -m+\epsilon$, followed by the limit $\epsilon\rightarrow 0$. As in the previous sections, terms of order $\epsilon^{-2}$ and $\epsilon^{-1}$ arise, and, as before, these must vanish. The first of these reduces to a limiting version of the case \eqref{Prud9} discussed above; the second reproduces \eqref{5F4b1} with $a=-m$.

The limit $m\rightarrow -1$ yields ,with the help of \eqref{5F4b1}, the known result
\begin{equation}
\mapleinline{inert}{2d}{sum((-gamma-Psi(k+1))/(k+1)^4, k = 0 .. infinity) = 1/6*Zeta(3)*Pi^2-2*Zeta(5)}{\[\displaystyle \sum _{k=0}^{\infty }{\frac {-\gamma-\Psi \left( k+1 \right) }{ \left( k+1 \right) ^{4}}}=\frac{1}{6}\,\zeta \left( 3 \right) {\pi}^{2}-2\,\zeta \left( 5 \right) \]} \,.
\label{Z1}
\end{equation}

Finally, reverting to the leading term (of order $\epsilon^0$) we initially find the reduction of an infinite series (on the left-hand side) to a collection of finite series on the right, all in terms of a single integer parameter $m$ (see \eqref{R3d} below where $m$ has been replaced by $a$ for reasons also discussed immediately below).\newline 

{\bf (Digression and Observation:} Experience has shown that occasionally, obscure underlying relationships exist that govern the behaviour of identities such as the following when $m$ is replaced by $a$. For a simple example of such a hidden association, see the discussion surrounding \eqref{Shpot1} - \eqref{Shpot3}, and for an extreme example, see Section \ref{sec:hail}. Consequently, a fastidious analyst should always check numerically for the possibility that another more  general identity may be lurking when a new identity involving $m \in \mathbb{N} $ is explored. The requirement is simply that the new identity be unambiguously calculable when subjected to the replacement $m:=a$, where $a \in \mathfrak{C}$. The goal is to discover if the more general identity might be numerically valid, all the while recognizing that there are an infinite number of ways this transition could be performed, none of which might be valid. The case where an identity is not unambiguously calculable is explored in Section \ref{sec:hail}. If the generalized identity appears to be valid, certification via the WZ algorithm \cite{WZ} becomes a possibility.) \newline

Keeping this possibility in mind, numerical tests of \eqref{R3d} below, obtained from the leading term (of order $\epsilon^0$) and derived under the condition that $a\rightarrow m \in \mathbb{N}$, indicate that it {\bf is} numerically true under the generalization $m:=a$. However, certification using the WZ algorithm \cite{WZ} fails due to the fact that both sides contain digamma functions, and hence are not of hypergeometric type. We are thus left with the following numerically tested, but analytically unproven identity (unless $a=m$). With reference to \eqref{D1S}, the general result ($c=-m+\epsilon$, of order $\epsilon^0$, $m:=a$) is:

\begin{align} \label{R3d}
&\displaystyle\sum _{k=0}^{\infty }{\frac {\Psi \left( 2+k \right) \Gamma \left( 1+k \right) }{ \left( 1+k \right) ^{2}\Gamma \left( k+3+a \right) }}=\frac{1}{2}\,\sum _{k=0}^{\infty }{\frac { \left( -1 \right) ^{k} \left( \Psi \left( 1+k \right) ^{2}-\Psi^{\prime} \left(1+k \right)  \right) }{ \left( 1+k \right) ^{3}\Gamma \left( 1+k \right) 
\mbox{}\Gamma \left( a-k+1 \right) }}
\\ \nonumber 
&+\frac {1}{\Gamma \left( 2+a \right) } \left(  \left( -  \Psi \left( 2+a \right)^{2}+\Psi^{\prime} \left( 2+a \right)  \right)  \left( 3/4\,{\gamma}^{2}-1/8\,{\pi}^{2} \right) 
\mbox{}-1/6\,\gamma\,  \Psi \left( 2+a \right) ^{3} \right. \\ \nonumber
& \left. + \left( 1/2\,\gamma\,\Psi^{\prime} \left(2+a \right) -{\gamma}^{3}+1/3\,{\pi}^{2}\gamma-2\,\zeta \left( 3 \right)  \right) \Psi \left( 2+a \right) \right. \\ \nonumber
& \left. -{\frac {5\,{\gamma}^{4}}{12}}+1/6\,{\pi}^{2}{\gamma}^{2}- \left( 1/6\,\Psi^{(2)} \left(2+a \right) +10/3\,\zeta \left( 3 \right)  \right) \gamma
\mbox{}+{\frac {5\,{\pi}^{4}}{144}} \right) 
\end{align}

As commented above, the series on the right-hand side is truncated at $m+1$ terms if $a=m$, and it can be identified as
\begin{equation}
\displaystyle \lim _{b\rightarrow 1}{\frac {\partial ^{2}}{\partial {b}^{2}}} \left( {\frac {{\mbox{$_5$F$_4$}(1,1,1,1,-a;\,2,2,2,b;\,1)}}{\Gamma \left( b \right) \Gamma \left( a+1 \right) }} \right) =-\sum _{k=0}^{\infty }{\frac { \left( -1 \right) ^{k}
\mbox{} \left( \Psi \left( k+1 \right)  ^{2}-\Psi^{\prime} \left(k+1 \right)  \right) }{\Gamma \left( k+1 \right)  \left( k+1 \right) ^{3}\Gamma \left( a-k+1 \right) }} \,.
\label{Rhs}
\end{equation} 
As well, the case $a=0$ yields the interesting Euler sum (also see \eqref{Z1} and \eqref{AnsSums})
\begin{equation}
\mapleinline{inert}{2d}{Sum(Psi(k+1)/((k+1)^3*(2+k)), k = 0 .. infinity) = -(gamma+1)*Zeta(3)-(1-(1/6)*Pi^2)*gamma+(1/360)*Pi^4+1}{\[\displaystyle \sum _{k=0}^{\infty }{\frac {\Psi \left( k+1 \right) }{ \left( k+1 \right) ^{3} \left( 2+k \right) }}=- \left( \gamma+1 \right) \zeta \left( 3 \right) 
\mbox{}- \left( 1-\frac{\pi^2}{6} \right) \gamma+{\frac {{\pi}^{4}}{360}}+1\]} \,.
\label{Z2}
\end{equation}

\subsubsection{Special cases: $d=-n$}

In the case $d=-n$, the $_3F_2(1)$ on the right-hand side of \eqref{Y1b2} is a limiting case of the known entry 25 of \cite{Milgram447} - that is (after evaluating the appropriate limits)
\begin{align} \label{C2A}
&\displaystyle {\mbox{$_3$F$_2$}(c-1,c-1,c-1;\,c,c+n+1;\,1)}=\frac {\pi\, \left( c-1 \right) 
\mbox{}\Gamma \left( c+n+1 \right) }{\sin \left( \pi\,c \right) \Gamma \left( 2+n \right) }\\ \nonumber
&\times \left( {\frac { \left( -1 \right) ^{n}}{\Gamma \left( c-3-n \right) }\sum _{k=0}^{n+1}{\frac { \left( -1 \right) ^{k}\Psi \left( k+1 \right) }{\Gamma \left( n-k+2 \right)  \left( c-k-2 \right) \Gamma \left( k+1 \right) 
\mbox{}}}}+{\frac {\Psi \left(c -1 \right) }{\Gamma \left( c-1 \right) }} \right) \,.
\end{align}

Substitution of \eqref{C2A} into \eqref{Y1b2} then gives the new result

\begin{align} \label{C2A1}
&\displaystyle {\mbox{$_4$F$_3$}(1,1,c,-n;\,2,2,2;\,1)}=-{\frac {\Gamma \left( c-1 \right)  \left( -1 \right) ^{n}}{ \left( c-1 \right) \Gamma \left( c-3-n \right)  \left( n+1 \right) }\sum _{k=0}^{n+1}{\frac { \left( -1 \right) ^{k}\Psi \left( k+1 \right) }{\Gamma \left( n-k+2 \right)  \left(c -k-2 \right) \Gamma \left( k+1 \right) 
\mbox{}}}} \\ \nonumber
&\hspace{4cm} \mbox{}+{\frac {(\Psi \left( 2+n \right)+2\,\gamma) }{ \left( c-1 \right)  \left( 1+n \right) }}
\end{align}
and, after differentiating with respect to $c$

\begin{align} \label{C2Ad}
&\displaystyle \sum _{k=0}^{n }{\frac {\Gamma \left( k+1 \right) \Psi \left( c+k \right) \Gamma \left( c+k \right)  \left( -1 \right) ^{k}}{\Gamma \left( n-k+1 \right) 
\Gamma \left( k+2 \right)  ^{3}}}=\frac { \left( -1 \right) ^{n} \left( \Psi \left( c-3-n \right)-2\,\Psi \left(c -1 \right)   \right) 
\mbox{} \Gamma \left(c -1 \right)  ^{2}}{\Gamma \left( 2+n \right) \Gamma \left( c-3-n \right) } \\ \nonumber
&\times\sum _{k=0}^{n+1}{\frac { \left( -1 \right) ^{k}\Psi \left( k+1 \right) }{\Gamma \left( n-k+2 \right)  \left(c -k-2 \right) \Gamma \left( k+1 \right) }} +{\frac { \left( \Psi \left( 2+n \right) +2\,\gamma \right) \Psi \left(c -1 \right) \Gamma \left(c -1 \right) }{\Gamma \left( 2+n \right) }} \\ \nonumber
&+{\frac { \Gamma \left(c -1 \right)  ^{2}
\mbox{} \left( -1 \right) ^{n}}{\Gamma \left( 2+n \right) \Gamma \left( c-3-n \right) }\sum _{k=0}^{n+1}{\frac { \left( -1 \right) ^{k}\Psi \left( k+1 \right) }{\Gamma \left( n-k+2 \right)  \left( c-k-2 \right) ^{2}\Gamma \left( k+1 \right) }}}\\
&
\end{align}
or, alternatively (by  solving and substituting)
\begin{align} \label{C2bd}
&\displaystyle {\mbox{$_4$F$_3$}(1,1,c,-n;\,2,2,2;\,1)}=\\ \nonumber
&{\frac {\Gamma \left( c-1 \right)  \left( -1 \right) ^{n}}{\Gamma \left( c-3-n \right)  \left(\Psi \left( c-3-n \right)  -2\,\Psi \left( c-1 \right) \right) 
\mbox{} \left( c-1 \right)  \left( n+1 \right) }\sum _{k=0}^{n+1}{\frac { \left( -1 \right) ^{k}\Psi \left( k+1 \right) }{\Gamma \left( n-k+2 \right)  \left( c-k-2 \right) ^{2}
\mbox{}\Gamma \left( k+1 \right) }}}\\ \nonumber
&-{\frac {\Gamma \left( 1+n \right) }{ \left( -2\,\Psi \left( -1+c \right) +\Psi \left( c-3-n \right)  \right) 
\mbox{}\Gamma \left(c \right) }\sum _{k=0}^{n}{\frac {\Gamma \left( k+1 \right) \Psi \left( c+k \right) \Gamma \left( c+k \right)  \left( -1 \right) ^{k}}{\Gamma \left( n-k+1 \right) \Gamma \left( k+2 \right)^{3}}}}
\mbox{}\\ \nonumber
&+{\frac { \left( \Psi \left( -1+c \right) -\Psi \left( c-3-n \right)  \right)  \left( \Psi \left( 2+n \right) +2\,\gamma \right) }{ \left( 2\,\Psi \left( c-1 \right) -\Psi \left( c-3-n \right)  \right) 
\mbox{} \left( c-1 \right)  \left( n+1 \right) }}\,.
\end{align}

\subsubsection{Special case d=-c-n}

In the case that $d=-c-n$, the $_3F_2(1)$ on the right-hand side of \eqref{Y1b2} is easily identified ({\bf Note:} here only, $n\geq -2$) as being contiguous to Watson's theorem (\cite[entry 1]{Milgram447}) for which general closed forms are well known \cite{ChuW}. However, because the top and bottom parameters are almost all separated by integers, this becomes a very complicated general limiting case; when the limits are all evaluated, the result contains far too many terms to realistically reproduce here. For the first few values corresponding to the variable $n= -1,0,...$, the results are
\begin{equation}
\mapleinline{inert}{2d}{Cnm1 := hypergeom([c-1, c-1, c-1], [c, 2*c-1], 1) = (((1/2)*Psi(1, c)-(1/4)*Psi(1, (1/2)*c))*(c-1)^2+1/4)*4^c*GAMMA(-1/2+c)/(sqrt(Pi)*GAMMA(c))}{\[\displaystyle{\mbox{$_3$F$_2$}(c-1,c-1,c-1;\,c,2\,c-1;\,1)}={\frac {
\mbox{}{4}^{c}\Gamma \left( -1/2+c \right) }{ \sqrt{\pi}\Gamma \left( c \right) }}\]} \left(  \left( \frac{1}{2}\Psi^{\prime} \left( c \right) -\frac{1}{4}\Psi^{\prime} \left(c/2 \right)  \right)  \left( c-1 \right) ^{2}+1/4
\mbox{} \right) 
\label{Cnm1}
\end{equation} 
\begin{align} \label{Cn0}
&\displaystyle {\mbox{$_3$F$_2$}(c-1,c-1,c-1;\,c,2\,c;\,1)}=\frac {{4}^{c}\Gamma \left( 1/2+c \right) }{ \sqrt{\pi}\Gamma \left( c+1 \right) } \left(  \left( \Psi^{\prime} \left( c \right) -1/2\,\Psi^{\prime} \left( c/2 \right)  \right)  \left( c-1 \right) ^{2}-{\frac {-{c}^{2}+c-1}{2\,{c}^{2}}}
\mbox{} \right)\\ \nonumber
&\displaystyle {\mbox{$_3$F$_2$}(c-1,c-1,c-1;\,c,2\,c+1;\,1)}= \\ \label{Cn1}
&\hspace{3cm}\frac {{4}^{c}\Gamma \left( 1/2+c \right) }{ \sqrt{\pi}
\mbox{}\Gamma \left( c+2 \right) } \left(  c\,\left( -\Psi^{\prime} \left(c/2 \right)+2\,\Psi^{\prime} \left( c \right) \right)  \left( c-1 \right) ^{2}
\mbox{}+{\frac {{c}^{4}+{c}^{3}-2\,{c}^{2}+3\,c+1}{c \left( c+1 \right) ^{2}}} \right) \\ \nonumber
&\displaystyle \mbox{$_3$F$_2$} \left( c-1,c-1,c-1;\,c,2\,c+2;1\, \right) = 
\frac {{4}^{c}\Gamma \left( 3/2+c \right) 
\mbox{}}{ \sqrt{\pi}\Gamma \left( c+3 \right) } \\ \label{Cn2}
&\times \left( -2\,c \left( \Psi^{\prime} \left(c/2 \right) -2\,\Psi^{\prime} \left( c \right)  \right)  \left( c-1 \right) ^{2}
\mbox{}-{\frac {-2\,{c}^{6}-10\,{c}^{5}-12\,{c}^{4}+14\,{c}^{3}-10\,{c}^{2}-44\,c-8}{ \left( c+2 \right) ^{2} \left( c+1 \right) ^{2}c}} \right)\,. 
\end{align}
The case $n=-2$ equates to \eqref{Entry} with $c=a+1$. The corresponding results for respective values of $n=-1,0,\dots$   applied to \eqref{Y1b2} are as follows:

\begin{align} \label{Y3m1}
\displaystyle{\mbox{$_4$F$_3$}(1,1,c,2-c;\,2,2,2;\,1)}=&{\frac {1}{ \left( c-1 \right) ^{2}} \left( 2\,\gamma+2\,\Psi \left( c \right) +{\frac { \left( \Psi^{\prime} \left(c/2 \right) -2\,\Psi^{\prime} \left(c \right)  \right) \sin \left( \pi\,c \right) }{\pi}}
\mbox{} \right) }\\ \nonumber
&- \frac{1}{\left( c-1 \right) ^{3}}-{\frac {\sin \left( \pi\,c \right) }{\pi\, \left( c-1 \right) ^{4}}}
\end{align}

\begin{align} \label{Y30}
\displaystyle {\mbox{$_4$F$_3$}(1,1,c,1-c;\,2,2,2;\,1)}=&{\frac {\sin \left( \pi\,c \right) }{\pi} \left( {\frac {\Psi^{\prime} \left( c/2 \right) -2\,\Psi^{\prime} \left( c \right) }{c \left( c-1 \right) }}-{\frac {{c}^{2}-c+1}{{c}^{3} \left( c-1 \right) ^{3}}} \right) }
\\ \nonumber
&+2\,{\frac {(\Psi \left( c+1 \right)+\gamma) }{c \left( c-1 \right) }}-{\frac {2\,c-1}{{c}^{2} \left( c-1 \right) ^{2}}}
\end{align}

\begin{align} \label{Y31}
&\displaystyle {\mbox{$_4$F$_3$}(1,1,c,-c;\,2,2,2;\,1)}=2\,{\frac {\Psi \left( c+2 \right)+\gamma }{ \left( c-1 \right)  \left( c+1 \right) }}+{\frac {1-3\,{c}^{2}}{c \left( c+1 \right) ^{2} \left( c-1 \right) ^{2}}}\\ \nonumber
& -\frac {\Gamma \left( c \right) \sin \left( \pi\,c \right) }{\pi\,\Gamma \left( c+2 \right)  \left( c-1 \right) ^{3}}
 \left(c \left( 2\,\Psi^{\prime} \left(c \right)-\Psi^{\prime} \left(c/2 \right) \right)  \left( c-1 \right) ^{2}
\mbox{}+{\frac {{c}^{4}+{c}^{3}-2\,{c}^{2}+3\,c+1}{c \left( c+1 \right) ^{2}}} \right) 
\end{align}

\begin{align} \label{Y32}
&\displaystyle {\mbox{$_4$F$_3$}(1,1,c,-1-c;\,2,2,2;\,1)}=2\,{\frac {\Psi \left( 3+c \right) +\gamma}{ \left( c-1 \right)  \left( c+2 \right) }}-2\,{\frac { \left( 2\,c+1 \right)  \left( {c}^{2}+c-1 \right) }{c \left( c+1 \right)  \left( c+2 \right) ^{2} \left( c-1 \right) ^{2}}} \\ \nonumber
&+\frac {\Gamma \left( c \right) \sin \left( \pi\,c \right) {4}^{c}\Gamma \left( 3/2+c \right) }{{\pi}^{3/2} \left( c+2 \right)  \left( c-1 \right) ^{3}
\mbox{}\Gamma \left( 2\,c+2 \right) } \left( 2\,c \left( \Psi^{\prime} \left(c/2 \right) -2\,\Psi^{\prime} \left(c \right)  \right)  \left( c-1 \right) ^{2}  \right. \\ \nonumber 
&\hspace{5cm}\left.+{\frac {-2\,{c}^{6}-10\,{c}^{5}-12\,{c}^{4}+14\,{c}^{3}-10\,{c}^{2}-44\,c-8}{ \left( c+2 \right) ^{2}c \left( c+1 \right) ^{2}}}
 \right) \,. 
\end{align}
The case $n=-2$ has been omitted here also, because it corresponds to \eqref{4F3} with $a=2-c$, providing an alternative direct derivation of \eqref{4F3}. 
\subsubsection{Special Case d=n-c}

In this special case, the series representation on the left-hand side of \eqref{Y1b2} will only converge if $n\leq3$. Since the case $n=3$ overlaps with the special cases listed in the previous subsection, for that case we find \eqref{Entry}, and correspondingly an equivalent to \eqref{4F3}.

\subsection{c=1}

Setting $c=1$ in \eqref{4F3} gives a known result on the left (see \eqref{Prud9}), and after a limiting calculation on the right, an interesting transformation results:


\begin{align} \label{Y2}
&\displaystyle {\mbox{$_4$F$_3$}(1,1,1,1;\,2,2,d;\,1)}={\frac {1}{\Gamma \left( 1-d
\mbox{} \right) }\sum _{k=0}^{\infty }{\frac {\Psi \left( 1+k \right) \Gamma \left( 3-d+k \right) }{ \left( 1+k \right) ^{3}\Gamma \left( k+1 \right) }}}+\frac{d-1}{6}\Psi \left( d \right)^{3} \\ \nonumber
&+ \left(\frac{d-1}{12} \left( -{\pi}^{2}+18\,{\gamma}^{2}-6\,\Psi^{\prime} \left(d \right)  \right) -2\,\gamma \right) \Psi \left( d \right) + \left( 1/2-\left(d-1 \right) \gamma \right)( \Psi^{\prime} \left(d \right)-\Psi(d)^2) 
\mbox{} \\ \nonumber
&+\frac{d-1}{6} \Psi^{(2)} \left(d \right)+2/3\, \left( d-1 \right) {\gamma}^{3}-3/2\,{\gamma}^{2}+4/3\,\zeta \left( 3 \right) d
\mbox{}+1/12\,{\pi}^{2}-4/3\,\zeta \left( 3 \right). 
\end{align}

Comparing the right-hand sides of \eqref{Y2} and \eqref{Eq2p2}, with the aid of \eqref{5F4b1} identifies the following sum

\begin{align} \label{NeqA}
\displaystyle \sum _{k=0}^{\infty }\frac {\Gamma \left( 2-d+k \right) }{\Gamma \left( 1+k \right)  \left( 1+k \right) ^{2}
\mbox{}}& \left( \Psi^{\prime} \left(1+k \right) -{\frac {\Psi \left( 1+k \right) }{1+k}} \right) =  
\\ \nonumber
& \Gamma \left( 1-d \right) \bigg{(} -{\Psi \left( d \right)^{3}}/3-3\Psi \left( d \right) ^{2}\gamma/2- \left( 2\,{\gamma}^{2}+{\pi}^{2}/6 \right) \Psi \left( d \right) 
+\Psi^{\prime} \left(d \right) \gamma/2  \\ \nonumber 
& +\Psi^{(2)} \left(d \right)/6 -5{\gamma}^{3}/6-{\pi}^{2}\gamma/4\,+\zeta \left( 3 \right)/3  \bigg{)}
\end{align}

In the case that $d=n$, replace the term containing the sum on the right-hand side of \eqref{Y2} with
\begin{equation}
\Gamma \left( n \right) \sum _{k=0}^{n-3}{\frac { \left( -1 \right) ^{k}\Psi \left( 1+k \right) }{ \left( 1+k \right) ^{3}\Gamma \left( 1+k \right) \Gamma \left( n-2-k \right) }}.
\label{Yna}
\end{equation}
All other terms are the same subject to the replacement $d=n$. \eqref{Y2} and \eqref{Eq2p2} are unusual, in that they represent one of the few transformations of a hypergeometric function containing a free bottom parameter.

\section{New Transformations} \label{sec:new}

By comparing several of the foregoing results, a number of new transformations were found. Specifically:
\begin{itemize}
\item Comparison of \eqref{Miller1Final} and \eqref{Sm3d} produces the identity

\begin{align} \label{S1Id}
&\displaystyle {\mbox{$_5$F$_4$}(1,1,1,2-a,a+1;\,2,2,2,2;\,1)}=-{\frac {\sin \left( \pi\,a \right) S_{1}(a)}{3\,a \left( a-1 \right) \pi}}
\\ \nonumber
&-{\frac {2}{3\,a \left( a-1 \right) } \left( {\gamma}^{2}+ \left( \Psi \left( a \right) +\Psi \left( 1-a \right)  \right) \gamma+ \left( \Psi \left( a \right)  \right) ^{2}+{\frac {\pi\,\Psi \left( a \right) \cos \left( \pi\,a \right) }{\sin \left( \pi\,a \right) }}
\mbox{}-\frac{{\pi}^{2}}{6} \right) }\,.
\end{align}

\item Comparison of \eqref{Miller1Final} with \eqref{Heq1} using \eqref{Entry} yields a related transformation

\begin{align} \label{5F4}
&\displaystyle {\mbox{$_5$F$_4$}(1,1,1,2-a,a+1;\,2,2,2,2;\,1)}={\frac { \left( \Psi^{\prime} \left( a/2+1/2 \right) -\Psi^{\prime} \left( a/2 \right)  \right) 
\mbox{}\sin \left( \pi\,a \right) }{2\,{a}^{2} \left( a-1 \right) \pi}}+2\,{\frac { \left( \gamma+\Psi \left( a \right)  \right) ^{2}}{a \left( a-1 \right) }} \\ \nonumber
&+ {\frac {{4}^{-a}\sin \left( \pi\,a \right) \Gamma \left( a \right) }{a \sqrt{\pi} \left( a+1 \right) ^{2}
\mbox{}\Gamma \left( a+1/2 \right)  \left( a-1 \right) }}
{\mbox{$_4$F$_3$}(a+1,a+1,a+1,a+1;\,2\,a+1,a+2,a+2;\,1)} \,.
\end{align}

\item A simple comparison between \eqref{S1Id} and \eqref{5F4} gives

\begin{align} \label{R12}
&\displaystyle {\mbox{$_4$F$_3$}(a+1,a+1,a+1,a+1;\,a+2,a+2,2\,a+1;\,1)}=-{\frac {\Gamma \left( a+1/2 \right) {4}^{a} \left( a+1 \right) ^{2}{ S_{1}(a)}}{ 3\,\sqrt{\pi}
\mbox{}\Gamma \left( a \right) }} \\ \nonumber
&+\frac {\Gamma \left( a+1/2 \right) {4}^{a} \left( a+1 \right) ^{2}}{3\,\Gamma \left( a \right) } \left( -8\,{\frac { \sqrt{\pi}{\gamma}^{2}}{\sin \left( \pi\,a \right) }}-{\frac {2\, \sqrt{\pi} \left( 7\,\Psi \left( a \right) +\Psi \left( 1-a \right)  \right) \gamma}{\sin \left( \pi\,a \right) }} \right. \\ \nonumber
&\left. -\,\frac{3}{2}{\frac {\left( \Psi^{\prime} \left(a/2+1/2 \right) -\Psi^{\prime} \left(a/2 \right)\right) }{a \sqrt{\pi}}}+{\frac {\displaystyle \frac{{\pi}^{5/2}}{3}-8\, \sqrt{\pi} \Psi \left( a \right)  ^{2}}{\sin \left( \pi\,a \right) }}-2\,{\frac {{\pi}^{3/2}\cos \left( \pi\,a \right) \Psi \left( a \right) }{ \sin^2 \left( \pi\,a \right)  }}
\mbox{} \right) \,.
\end{align}

\item Comparing \eqref{Eq4} and \eqref{Sm3d} with $a<1$ gives a relationship between two sums:

\begin{align} \label{SxS1}
&\displaystyle \sum _{k=0}^{\infty }{\frac {\Gamma \left( 1+k \right) \Psi \left( 1-a+k \right) }{\Gamma \left( 1-a+k \right)  \left( 1-a+k \right)^2 }}={\frac {\sin \left( \pi\,a \right) \Gamma \left( a \right) { S_1(a)}}{3\,\pi}}
\mbox{}+\frac{2\,\Psi \left( a \right) \Gamma \left( a \right) \gamma}{3}+\frac{14}{3}\Psi \left( 1-a \right) \Gamma \left( a \right) \gamma\\ \nonumber
& 
+\frac{8}{3}\Gamma \left( a \right) {\gamma}^{2}+ \left( {\frac {\Psi \left( 1-a \right)  \left( \Psi^{\prime} \left( a/2 \right) -2\,\Psi^{\prime} \left(a \right)  \right) \sin \left( \pi\,a \right) }{\pi}}-\frac{1}{9}{\pi}^{2}+\frac{8}{3} \Psi \left( a \right)\,\Psi \left( 1-a \right)  
 \right) \Gamma \left( a \right) \,.
\end{align}

\item After operating on \eqref{SxS1} with $a \rightarrow 1-a$, and noting that $S_1(a)$ is invariant under that operation, eliminate $S_1(a)$ from the two resulting equations, and obtain a transformation between two very slowly converging sums:


\begin{align} \label{Sx3}
&\displaystyle \Gamma \left( a \right) \sum _{k=0}^{\infty }{\frac {\Gamma \left( 1+k \right) \Psi \left( a+k \right) }{\Gamma \left( a+k \right)  \left( a+k \right) ^{2}}}=\Gamma \left( 1-a \right) \sum _{k=0}^{\infty }{\frac {\Gamma \left( 1+k \right) \Psi \left( 1-a+k \right) }{\Gamma \left( 1-a+k \right)  \left( 1-a+k \right) ^{2}}}-{\frac {{4\,\pi}^{2}\cos \left( \pi\,a \right) \gamma}{  \sin^2 \left( \pi\,a \right)  }} \\ \nonumber
&+ \left( \Psi^{\prime} \left( 1/2-a/2 \right) -2\,\Psi^{\prime} \left(1-a \right)  \right) \Psi \left( a \right) + \left( 2\,\Psi^{\prime} \left(a \right) -\Psi^{\prime} \left(a/2 \right)  \right) \Psi \left( 1-a \right) \,.
\end{align}

The left-hand side of \eqref{Sx3} converges for $a>0$ while the sum on the right-hand side converges for $a<1$; since an overlap exists, the two sides of the equation are analytic continuations of each other. That is, the right-hand side represents the left hand side for $a\leq0$ and, by reversing the equation, the reverse will be true for $a\geq 1$. Both series are very slow to converge when $0<a<1$, and it is very difficult to verify \eqref{Sx3} numerically in the overlap range of the variable $a$, except at $a=1/2$, when the two sides reduce to an identity.\newline

\item After operating on both sides of \eqref{4F3} with $\frac{\partial}{\partial\,a}$ the following sum is found


\begin{align} \label{S1Sum}
&\displaystyle \sum _{k=1}^{\infty }{\frac {\Gamma \left( 1-a+k \right) \Gamma \left( a+k \right)  \left( \Psi \left( a+k \right) -\Psi \left( 1-a+k \right)  \right) }{ \Gamma \left(k \right)^{2}
\,k^3}} \\ \nonumber
&\hspace{2cm}=2\,\Psi^{(2)} \left(a \right) -\frac{1}{2}\Psi^{(2)}\, \left(a/2 \right) -{\frac {2\,\pi\,\Psi^{\prime} \left(a \right) }{\sin \left( \pi\,a \right) }}
+{\frac {2\, \left( \gamma+\Psi \left( a \right)  \right) {\pi}^{2}\cos \left( \pi\,a \right) }{ \sin^2 \left( \pi\,a \right)  }}
\end{align}

A comparison of \eqref{S1Sum} with \eqref{S1} suggests a redefinition of the sum $S_1(a)$:
\begin{align} \label{S1_new}
\displaystyle {S_{1}}(a)&=-2\,\sum _{k=1}^{\infty }{\frac {\Gamma \left( k+1-a \right) \Gamma \left( a+k \right)  \left( \Psi \left( k \right) -\Psi \left( a+k \right)  \right) }{{k}^{3} \Gamma \left( k \right) ^{2}
\mbox{}}}+\frac{1}{2}\Psi^{(2)} \left( a/2 \right) -2\,\Psi^{(2)} \left( a \right)\\ \nonumber
& +{\frac {2\,\pi\,\Psi^{\prime} \left(a \right) }{\sin \left( \pi\,a \right) }}-{\frac { 2\,\left( \gamma+\Psi \left( a \right)  \right) {\pi}^{2}\cos \left( \pi\,a \right) }{ \sin^2 \left( \pi\,a \right) }}
\end{align}

\item
To aid in the verification of \eqref{AnsIm}, we explore the reduction of that result in the limit $a\rightarrow n$. This requires that the sums be split at $k=n-1$ and appropriate limits be evaluated. The calculation is straightforward, and requires the invocation of \eqref{Lemma4} and \eqref{Prud9} with $a=1$ and $b=2-n$, all of which eventually give the  following representation reducing an infinite to finite sums:  
\begin{align} \label{AnsImG}
&\displaystyle \sum _{k=0}^{\infty }{\frac {\Gamma \left( 1+k \right) \Psi \left( k+n \right) }{ \left( k+n \right) ^{2}\Gamma \left( k+n \right) }}\\ \nonumber 
&=\Psi \left( n \right) \sum _{k=1}^{n-1}{\frac {\Gamma \left( n-k \right) \Psi \left( k \right)  \left( -1 \right) ^{k}}{  \Gamma \left( n-k+1 \right) ^{2}\Gamma \left( k \right) }}
\mbox{}+ \left( -1 \right) ^{n}\sum _{k=1}^{n-1}{\frac { \left( -1 \right) ^{k}\Psi \left( n-k \right)  \left( \Psi \left( k \right) -\Psi \left( n \right)  \right) }{{k}^{2}\Gamma \left( k \right) \Gamma \left( n-k \right) }}\\ \nonumber
&- \frac{\left( -1 \right) ^{n}}{\Gamma \left( n \right)} { { \left( \Psi \left( n \right)^{3}/2+2\,\gamma\, \Psi \left( n \right)^{2}- \left({\pi}^{2}/12\,-3{\gamma}^{2}/2+1/2\,\Psi^{\prime} \left(n \right)  \right) \Psi \left( n \right) 
-2\,\gamma\,\Psi^{\prime} \left(n \right) -1/2\Psi^{(2)} \left(n \right)  \right) }{}}\,.
\end{align}
Setting $n=1\dots4$ gives the first few specific results for a family of generalized Euler sums
\begin{subequations}
\begin{align} \label{AnsSums} \nonumber
& \\ 
&\hspace{2cm}\displaystyle \sum _{k=0}^{\infty }{\frac {\Psi \left( 1+k \right) }{ \left( 1+k \right) ^{2}}}=-{\pi}^{2}\gamma/6+\zeta \left( 3 \right) \\  
&\hspace{1.5cm}\displaystyle \sum _{k=0}^{\infty }{\frac {\Psi \left( 2+k \right) }{ \left( 1+k \right)  \left( 2+k \right) ^{2}}}= \left( 1+\gamma \right){\pi}^{2}/6 
\mbox{}-2\,\gamma-\zeta \left( 3 \right)    \\
&\hspace{1cm}\displaystyle \sum _{k=0}^{\infty }{\frac {\Psi \left( k+3 \right) }{ \left( 1+k \right)  \left( 2+k \right)  \left( k+3 \right) ^{2}}}= \left( -1/8-\gamma/12 \right) {\pi}^{2}
\mbox{}+3/4\,\gamma+3/4+\zeta \left( 3 \right)/2   \\
&\displaystyle \sum _{k=0}^{\infty }{\frac {\Psi \left( k+4 \right) }{ \left( 1+k \right)  \left( 2+k \right)  \left( k+3 \right)  \left( k+4 \right) ^{2}}}= \left( {\frac{11}{216}}+\gamma/36 \right) {\pi}^{2}
\mbox{}-{\frac {31\,\gamma}{108}}-{\frac{5}{18}}
\mbox{}-\zeta \left( 3 \right)/6 ,
\end{align}
\end{subequations}
the first of which is well known. See also \eqref{SxS1} and \eqref{Lemma2}.

\item
Although \eqref{Lemma4} and \eqref{ht} were quoted independently and intended for different uses, in the case that $a=n$, \eqref{ht} reduces to \eqref{Lemma4}, in which case a comparison of the right-hand sides (with help from \eqref{Prud9} using $a=1$ and $b=2-n$),  eventually yields
\begin{align} \label{Newsum2}
\displaystyle \sum _{k=0}^{n}&{\frac { \left( -1 \right) ^{k}\Psi \left( 1+k \right) }{\Gamma \left( n-k+1 \right)  \left( n-k+1 \right) ^{2}\Gamma \left( 1+k \right) 
\mbox{}}}={\frac { \left(  \Psi \left( n+2 \right) ^{2}+\gamma\,\Psi \left( n+2 \right) -\Psi^{\prime} \left(n+2 \right)  \right)  \left( -1 \right) ^{n}
\mbox{}}{\Gamma \left( n+2 \right) }}\\ \nonumber
&\hspace{1cm}-{\frac {\Psi^{\prime} \left(n/2+1 \right) -\Psi^{\prime} \left(n/2+3/2 \right) }{2\,\Gamma \left( n+2 \right) }}.
\end{align}

\item
The result \eqref{Sm1b} contained an arbitrary parameter $d$; in the case $d\rightarrow 0$, as has been seen here frequently, the divergent term of order $1/d$ in the expansion must vanish. This leads to


\begin{align} \label{Smm1H}
\displaystyle\mbox{$_4$F$_3$}(a,1-a,a-b,a-b;&\,1,1+a-b,1+a-b;\,1)=\\ \nonumber
&{\frac { \sin \left( \pi\,a \right)  \left( a-b \right) ^{2}}{\pi}}S_{{2}} \left( a,b \right)
+{\frac {\pi\, \Gamma \left( 1+a-b \right) ^{2}}{\Gamma \left( 1-b \right) \sin \left( \pi\,a \right) \Gamma \left( 2\,a-b \right) }}
\end{align}

In the above, if $b\rightarrow a-1$ the left-hand side is well known (e.g. \cite[Entry 26]{Milgram447}) leading to the identification of the sum
\begin{align} \label{Smm1A}
\displaystyle \sum _{k=1}^{\infty }&{\frac {\Gamma \left( a+k-1 \right) \Gamma \left(k -a \right)  \left( \Psi \left( k-a \right) -2\,\Psi \left( k \right) +\Psi \left( a+k-1 \right)  \right) 
\mbox{}}{\Gamma \left( k \right)^{2}k}}\\ \nonumber
&\hspace{2cm}={\frac {1}{ \left( a-1 \right) ^{2}{a}^{2}}}+{\frac {2\,\pi}{a\sin \left( \pi\,a \right)  \left( -1+a \right) }}
\end{align}
The  case $b\rightarrow a$ reproduces \eqref{S1Id}.

\item 
Differentiate \eqref{Smm1H} with respect to $b$ then evaluate the limit $b\rightarrow a$ to produce the identity


\begin{align} \label{6F5}
&\displaystyle a \left( a-1 \right) {\mbox{$_6$F$_5$}(1,1,1,1,2-a,1+a;\,2,2,2,2,2;\,1)}=-{\frac {\sin \left( \pi\,a \right) }{4\,\pi}S_{{4}}(a)}-\frac{\zeta(3)}{6} -\frac{\Psi^{(2)}\!\left(a \right)}{12}\\ \nonumber
&-{\frac {\pi\, \left( \gamma+\Psi \left( a \right)  \right) ^{2}\cos \left( \pi\,a \right) }{2\,\sin \left( \pi\,a \right) }}
+\frac{1}{6} \left( \gamma+\Psi \left( a \right)  \right) {\pi}^{2}-\frac{1}{3} \left( \gamma+\Psi \left( a \right)  \right) ^{3}
\end{align}

where
\begin{equation}
\mapleinline{inert}{2d}{S[2] = Sum(GAMMA(a+k)*GAMMA(k+1-a)*(Psi(k+1-a)+Psi(a+k)-2*Psi(k+1))/(k^4*GAMMA(k)^2), k = 1 .. infinity)}{\[\displaystyle S_{{4}}(a)=\sum _{k=1}^{\infty }{\frac {\Gamma \left( a+k \right) \Gamma \left( k+1-a \right)  \left( \Psi \left( k+1-a \right) +\Psi \left( a+k \right) -2\,\Psi \left( k \right)  \right) 
\mbox{}}{{k}^{4}\, \Gamma \left( k \right) ^{2}}}\]}\,.
\label{S4}
\end{equation}

Alternatively, let $b\rightarrow a-1$, and, with reference to \eqref{Lemma3b}, obtain a closed expression for the infinite sum

\begin{align} \label{Smm1dL0}
&\displaystyle \sum _{k=1}^{\infty }{\frac { \left( \Psi \left( k+1-a \right) -2\,\Psi \left( k+1 \right) +\Psi \left( a+k \right)  \right) \Gamma \left( a+k \right) 
\mbox{}\Gamma \left( k+1-a \right) }{ \Gamma \left( k+2 \right)  ^{2}}}\\ \nonumber
&=2\,{\frac {\Psi^{\prime} \left(a/2 \right) }{a \left( a-1 \right) }}
-2\,{\frac {\pi\, \left( {a}^{2}-a-3 \right) \Psi \left( a \right) }{a \left( a-1 \right) \sin \left( \pi\,a \right) }}
\mbox{}-4\,{\frac {\Psi^{\prime} \left(a \right) }{a \left( a-1 \right) }}-{\frac {{\pi}^{2} \left( {a}^{2}-a-1 \right) \cos \left( \pi\,a \right) }{a \left( a-1 \right)  \sin^2 \left( \pi\,a \right) }}\\ \nonumber
&-{\frac {2\,{a}^{2}-2\,a+2}{{a}^{3} \left( a-1 \right) ^{3}}}
-{\frac {1}{\sin \left( \pi\,a \right) } \left( {\frac {2\,\pi\, \left( {a}^{2}-a-3 \right) \gamma}{a \left( a-1 \right) }}+{\frac {3\,\pi}{{a}^{2} \left( a-1 \right) ^{2}}} \right) }\,.
\end{align}
Cases with $b\rightarrow a-n$ are accessible by similar means.


\item Finally, comparing \eqref{OldEq8} and \eqref{Cx0d} gives the reduction of an infinite series to a finite sum:
\begin{align} \label{Eq82}
&\displaystyle {\frac {\Gamma \left( 2\,n+1 \right) {\mbox{$_5$F$_4$}(1,1,n,n,2\,n;\,n+1,n+1,n+1,n+1;\,1)}}{{n}^{5}  \Gamma \left( n \right)^{2}}}= \\ \nonumber
&\left( 2\,\Gamma \left( n \right) \sum _{k=1}^{n-1}{\frac { \left( -1 \right) ^{k}\Psi \left( n-k \right)  \left( \Psi \left( k \right) -\Psi \left( n \right)  \right) }{{k}^{2}\Gamma \left( n-k \right) \Gamma \left( k \right) }}+2\,\sum _{k=1}^{n-1}{\frac {\Gamma \left( k+n \right)  \left( -1 \right) ^{k} \left( \Psi \left( n-k \right) -\Psi \left( k+n \right)  \right) }{ \Gamma \left( k \right) ^{2}\Gamma \left( n-k \right) {k}^{4}}} \right. \\ \nonumber
&\left. - \left( 5\,\Psi \left( n \right) +4\,\gamma \right) \Psi^{\prime} \left( n \right) -3\,\Psi \left( n \right) ^{3}-6\, \Psi \left( n \right)  ^{2}\gamma+ \left( -3\,{\gamma}^{2}+{\pi}^{2}/6 \right) \Psi \left( n \right) 
\mbox{}+\Psi^{(2)} \left(n \right)  \right)  \left( -1 \right) ^{n}
\end{align}

\end{itemize}

\section{Are there any more?}

Apparently not! Although the previous two Sections yielded a rich assortment of identities, based on the operation \eqref{reduction} applied to \eqref{Miller1}, itself inspired by the discovery of \eqref{4F3}, when \eqref{reduction} was applied to either of the other 4-part transformations - \eqref{Miller2} or \eqref{GF_Imag} - only a trivial identity was found in both cases.\newline

\section{Summary}

This work is an attempt to summarize a number of techniques that are available to attempt to answer questions such as that posed by \eqref{question}. Although they are well known in the classical literature, three and four part relations are not usually applied to such problems. Typically the attempt produces relationships between various functions none of which can be evaluated in a useful form. In the cases studied here the utility of such relations devolves from the fact that difficult limits were required to be computed reducing the multi-part relations to simpler forms.\newline

Typically, the computation of such limits involves long and arduous calculations which occasionally yield a result that can be evaluated in a desirable form, as has been demonstrated - although not always the result sought! With the availability of modern computer tools and methods, many of these otherwise lengthy results, which heretofore have appeared in the literature as theorems in their own right, (e.g. \cite{SaiSax}), can now be obtained as needed. It is anticipated that the methods outlined here will be applicable to similar and more difficult problems. It is the author's opinion that the choice of the bottom parameter $2a$ in \eqref{question} is irrelevant to the specificity of the results obtained here. This suggests that more general closed sums of the form found in Sections 7 and 8 can be obtained by employing the methods outlined here using a generalized bottom parameter.\newline

It is interesting to note that the function $_2F_1(a,a;2a;1)$ is known by virtue of Gauss' theorem and the more generalized form $_3F_2(a,a,a;2a,a+1;1)$ is also known by virtue of Watson's theorem - see \eqref{Entry}. As demonstrated here, no further generalization is known, although a solution to $_3F_2(a,a,b;2a,b+1;1)$ would clearly provide a solution to \eqref{question}. In the closing statement of their paper \cite{Miller&Paris} Miller and Paris write: ``...there remains the open problem of deducing a summation formula for the series $_3F_2(a,b,f;c,f+n;1)$... We hope that the developments presented herein will stimulate further interest in this problem". Clearly the answer to the question discussed in this work would benefit from a solution to a limited subset of that problem. \eqref{Shpot4} may be a step in that direction.

\section{Acknowledgements}

The original question \eqref{question} was posed by Dalimil Mazac of the Simons Center for Geometry and Physics and Yang Institute for Theoretical Physics, Stony Brook, USA. I thank Larry Glasser, Mykola Shpot and anonymous referee(s) for comments on the original manuscript. 

\bibliographystyle{unsrt}

\bibliography{biblio}


\begin{appendices}
\section {4-part Transformations of a ${_4F_3(1)}$}
\subsection{Miller's transforms}
Miller \cite[Eqs. (1.1) and (1.2)] {Miller} has obtained the following two 4-part transformations among arbitrary ${_4F_3(1)}$. 

\begin{align} \label{Miller1}  
 {_4F_3}&(a,b,c,d;\,e,f,g;\,1)=  \Gamma \left( e \right) \Gamma \left( f \right) \Gamma \left( g \right) \Gamma \left( 1-d\right) \\ \nonumber
&\times  \left( {\frac {\Gamma \left( b-a \right) \Gamma \left( c-a \right) {_4F_3(a,1+a-e,1+a-f,1+a-g;\,1+a-d,1+a-b,1+a-c;\,1)}}{\Gamma \left( b \right) \Gamma \left( c \right) \Gamma \left( e-a \right)\Gamma \left( f-a \right) \Gamma \left( g-a \right) \Gamma \left( 1+a-d \right) }} \right. \\ \nonumber
& \left. +{\frac {\Gamma \left( a-b \right) \Gamma \left( c-b \right) {_4F_3(b,1+b-e,1+b-f,1+b-g;\,1+b-d,1+b-a,1+b-c;\,1)}}{\Gamma \left( a \right) \Gamma \left( c \right) \Gamma \left( e-b \right) \Gamma \left( f-b \right)\Gamma \left( g-b \right) \Gamma \left( 1+b-d \right) }} \right. \\ \nonumber
& \left.  +{\frac {\Gamma \left( a-c \right) \Gamma \left( b-c \right) {_4F_3(c,1+c-e,1+c-f,1+c-g;\,1+c-d,1+c-a,1+c-b;\,1)}}{\Gamma \left( b \right) \Gamma \left( a \right) \Gamma \left( e-c \right) \Gamma \left( f-c \right) \Gamma \left( g-c \right) \Gamma \left( 1+c-d \right) }} \right) 
\end{align}
and
\begin{align} \label{Miller2} 
&_4F_3(a,b,c,d;\,e,f,g;\,1)=\frac {\Gamma \left( e \right) \Gamma \left( f \right) \Gamma \left( 1-d \right) \Gamma \left( 1-c \right) }{\Gamma \left( 1-g \right) } 
\\ \nonumber
& \times
\left( { \frac {\Gamma \left( a-b \right) \Gamma \left( 1+b-g \right) {_4F_3(b,1+b-g,1+b-f,1+b-e;\,1+b-c,1+b-d,1+b-a;\,1)}}{\Gamma \left( a \right) \Gamma \left( e-b \right) \Gamma \left( f-b \right) \Gamma \left( 1+b-c \right) \Gamma \left( 1+b-d \right) }} \right. \\ \nonumber \left. \right. \\ \nonumber 
& 
\left. +  {\displaystyle \frac {\Gamma \left( b-a \right) \Gamma \left( 1+a-g \right) {_4F_3(a,1+a-g,1+a-e,1+a-f;\,1+a-c,1+a-d,1+a-b;\,1)}}{\Gamma \left( b \right) \Gamma \left( e-a \right) \Gamma \left( f-a \right) \mbox{}\Gamma \left( 1+a-c \right) \Gamma \left( 1+a-d \right) }} \right. \\ \left. \right.   \nonumber \\ 
&
\left. -  \frac {\Gamma \left(\! g-\!1 \right) \Gamma \left(\! 1\!+\!a\!-\!g \right) \Gamma \left(\! 1\!+\!b\!-\!g \right) {_{4}F_{3}}(1\!+\!b\!-\!g,\!1\!+\!a\!-\!g,\!1\!+\!c\!-\!g,\!1\!+\!d\!-\!g;\,2\!-\!g,1\!+\!e\!-\!g,\!1\!+\!f\!-\!g;\!1)}
{\Gamma \left( a \right) \Gamma \left( b \right) \Gamma \left( g-c \right) \Gamma \left( g-d \right) \Gamma \left( 1+e-g \right) \Gamma \left( 1+f-g \right) } \right) \nonumber
\end{align}
\subsection{Meijer G-function transforms} \label{sec:MeijerG}

Any hypergeometric function $_pF_q(a_p;b_q;z)$ can always be written as a Meijer G-function \cite[Chapter V]{Luke} and the reverse. The original intent of the Meijer G-function was to assign meaning to ${_pF_q(z)}$ when $p>q+1$ by the magic of contour integration and the residue theorem. This is summarized by the general result \cite[Eq. 5.3(1)] {Luke}
\begin{equation}
{\displaystyle  G^{m,n}_{p,q}\left(z\, \Big\vert\, {^{\displaystyle a_p}_{\displaystyle b_q}}\right)}=
                G^{n,m}_{q,p}\left(\frac{1}{z}\, \Big\vert\,^{\displaystyle 1-b_q}_{\displaystyle 1-a_p}\right) \hspace{1cm} arg\;\frac{1}{z}=-arg\;z\;.
\label{Ginv}
\end{equation}
which shows that the interchange of $a_p\leftrightarrow 1-a_p$ and $b_q\leftrightarrow 1-b_q$ together with $p\leftrightarrow q$ is equivalent to the analytic continuation $z\rightarrow 1/z$ for any ($m$ and $n$ combination(s) of) $_pF_{q-1}(z)$. Of course when $p=q$ and $z=1$ this is equivalent to a transformation among (combinations of) $_{q}F_{q-1}(1)$. For example, with $z=1$ in \eqref{Ginv} along with the prescription \cite[Eq. 5.2(7)]{Luke} for expanding a G-function into a combination of $_pF_{q-1}(1)$, and using $(-1)^a=\exp(-i\pi a)$ as a template for dealing with terms that arise in the expansion, rewrite both sides in the  form of (a) combination(s) of $_4F_3(1)$ by setting $m=1,n=4,p=q=4$. The imaginary part of the resulting expression must vanish, because for the moment, without loss of generality, it is  permissible to require that the variables $a, \dots g\in \mathbb{R}$. Solving that equation yields an apparently hitherto unrecognized 4-part transformation as follows  
\begin{align} \label{GF_Imag}
& {_4F_3} \left( a,b,c,d\,;e,f,g;\,1 \right) =-\frac {\Gamma \left( 1-a \right) 
\Gamma \left( 1-b \right) \Gamma \left( 1-d \right) \Gamma \left( 1-c \right) }{\Gamma \left( 1-f \right) \Gamma \left( 1-e \right) \Gamma \left( 1-g \right) }\times \\ \nonumber
&
\left(\!{\frac {\Gamma \left(g\! -1 \!\right) \Gamma \left( g-f \right) \Gamma \left(\! g-e \right) {_4F_3} \left(\! 1\!-g\!+a,1\!-g+b,1\!-g+c,1\!-g\!+\!d;2\!-g,f+1\!-g,e+1\!-g;\,1 \right)
}{\Gamma \left( g-b \right) \Gamma \left( g-c \right) \Gamma \left( g-d \right) \Gamma \left( g-a \right) }} \right. \\ \nonumber
& \left.
+{\frac {\Gamma \left(e \!-1 \right) \Gamma \left(e\! -f\! \right) \Gamma \left(e\! -g \right) {_4F_3} \left( \!1\!-e+\!a,1\!-e+b,1\!-e+c,1\!-e\!+d;2-e,f+1\!-e,\!g+1\!-e;\,1 \right)  
}{\Gamma \left( e-b \right) \Gamma \left( e-c \right) \Gamma \left( e-d \right) \Gamma \left( e-a \right) }} \right. \\ 
& \left.
+{\frac {\Gamma \left(f\! -1 \right) \Gamma \left( \!f\!-e \right) \Gamma \left( f\!-g \right) {_4F_3}\! \left( \!1\!-\!f\!+a,1\!-f\!+\!b,1\!-f+c,1\!-f\!+d;2\!-\!f,e\!+1\!-f,g\!+\!\!1\!-f;\,1\! \right) 
}{\Gamma \left( f-b \right) \Gamma \left( f-c \right) \Gamma \left( f-d \right) \Gamma \left( f-a \right) }} \right) \nonumber
\end{align}
whose validity can be expanded to include the variables $a,\dots g\ \in \mathbb{C}$ because of the principle of analytic continuation (parameter space is analytic). This relationship appears to be independent of \eqref{Miller1} and \eqref{Miller2}. However, it is easily shown that the transformation arising from the real part of that same expression is a symmetric permutation of, and  therefore equavalent to, \eqref{Miller1}. It is an open question which of the many other permutations of variables among \eqref{Miller1}, \eqref{Miller2} and \eqref{GF_Imag} will result in new, independent transformations applicable to $_4F_3(1)$, analogous to Whipple's categorization of three-part transformations for hypergeometric $_3F_2(1)$ discussed in \cite[Section 3.13]{Luke}. See Appendix C. This is discussed at some length in \cite{Raynal} and \cite{KrattRao}. \newline

\section{Application to a $_3F_2(1)$}

In Section \ref{sec:Evol}, the motivation for considering transformations among $_3F_2(1)$ was presented. With this in mind, several relationships involving special cases of $_3F_2(1)$ are (re)derived and collected in this Appendix.
 
\subsection{Special cases of Miller, Paris, Shpot and Srivastava} \label{sec:APPB}

A simple relevant transformation can be obtained using Miller and Paris \cite[Eq. (1.7)]{Miller&Paris} where it is shown for $m>n$ that
\begin{align} \label{Eq1p7a}
&{_3F_2} \left( a,b,n;\,c,m;\,1 \right) =\\ \nonumber
&{\frac {\Gamma \left( m \right) \Gamma \left( 1-a \right) \Gamma \left( 1-b \right) }{\Gamma \left( n \right) \Gamma \left( 1-c \right) }\sum _{k=0}^{m-n-1}{\frac { \left( -1 \right) ^{k}\Gamma \left( n+k \right) \Gamma \left( 1-c+n+k \right) }{\Gamma \left( k+1 \right) \Gamma \left( m-n-k \right) \Gamma \left( 1-a+n+k \right) 
\mbox{}\Gamma \left( 1-b+n+k \right) }}}
\mbox{}+ \\ 
&{\frac {\Gamma \left( m \right) \Gamma \left( c \right) \Gamma \left( 1-a \right) \Gamma \left( 1-b \right) }{\Gamma \left( c-a \right) \Gamma \left( c-b \right) \Gamma \left( m-n \right) }\sum _{k=0}^{n-1}{\frac { \left( -1 \right) ^{k}\Gamma \left( k+m-n \right) \Gamma \left( c-a-b+k+m-n \right) }{\Gamma \left( k+1 \right) \Gamma \left( n-k \right) \Gamma \left( 1-a+k+m-n \right) 
\mbox{}\Gamma \left( 1-b+k+m-n \right) }}}. \nonumber
\end{align}


Formally, generalize \eqref{Eq1p7a} by extending the upper limit of both sums to infinity (they are truncated by the denominator Gamma functions if $m,n \in \mathbb{N}$) straightforwardly identifying each of the extended sums as a hypergeometric function. Then relax the above condition by generalizing such that $m, n \in \mathbb{R}$, to obtain 
\begin{align} \label{Eq1p7b} 
&{_3F_2} \left( a,b,c;e,f;\,1 \right) =\\ \nonumber
&{\frac {\Gamma \left( f \right) \Gamma \left( 1-a \right) \Gamma \left( 1-b \right) \Gamma \left( -e+c+1 \right) {_3F_2} \left( c,-f+c+1,-e+c+1;1-b+c,1-a+c;\,1 \right) 
\mbox{}}{\Gamma \left( 1-e \right) \Gamma \left( f-c \right) \Gamma \left( 1-a+c \right) \Gamma \left( 1-b+c \right) }}
\\ \nonumber
&+{\frac {\Gamma \left( f \right) \Gamma \left( e \right) \Gamma \left( 1-a \right) \Gamma \left( 1-b \right) \Gamma \left( e-a-b+f-c \right) 
}{\Gamma \left( e-a \right) \Gamma \left( e-b \right) \Gamma \left( c \right) \Gamma \left( 1-a+f-c \right) \Gamma \left( 1-b+f-c \right) }}\\ 
&\times{_3F_2} \left( f-c,1-c,e-a-b+f-c;1-b+f-c,1-a+f-c;\,1 \right)  \,, \nonumber
\end{align}
a known result (see \eqref{P_and_M} below) that resides in \cite[Tables (3.3), (3.4) and(3.5)]{Luke}, valid for all values of the parameters $a,b,c,e,f$.  This observation justifies the apparently {\it ad hoc} replacement $c\rightarrow e, m\rightarrow f, n\rightarrow c$ in \eqref{Eq1p7a} to yield \eqref{Eq1p7b}. ({\bf Remark:} reverse this procedure by setting $f\rightarrow m$ and $c\rightarrow n$ to obtain a much simpler derivation of \eqref{Eq1p7a} from the known result \eqref{Eq1p7b}).\newline 

Similarly, Shpot and Srivastava \cite{Shpot&Sriv} have obtained a result for a more general related problem as follows:
\begin{align} \label{Shpot1}
&\displaystyle {\mbox{$_3$F$_2$}(a,b,c;\,b+1+m,c+1+n;\,z)}=\frac {\Gamma \left( b+1+m \right) \Gamma \left( c+1+n \right) }{  \Gamma \left( b \right) \Gamma \left( c \right) } \\ \nonumber
&\times \sum _{i=0}^{m} \sum _{j=0}^{n}\frac { \left( -1 \right) ^{i}
\mbox{}  \left( -1 \right) ^{j} }{\Gamma \left( 1-i+m \right) \Gamma \left( n+1-j \right) \Gamma \left( i+1 \right) 
\mbox{}\Gamma \left( j+1 \right) }  \\  \nonumber
&\times \left( {\frac {{\mbox{$_2$F$_1$}(a,c+j;\,c+j+1;\,z)}}{ \left( c+j \right)  \left( b-c+i-j \right) }}+{\frac {{\mbox{$_2$F$ _1$}(a,b+i;\,b+i+1;\,z)}}{ \left( b+i \right)  \left( c-b+j-i \right) }} \right)   
\end{align}

After a somewhat lengthy analysis, in the case that $z=1$, they reduce this result to a 3-part transformation \cite[Eq. (31)] {Shpot&Sriv} (see \eqref{SSEq31}) among terminating $_3F_2(1)$ on the right, and the same $_3F_2(1)$ that appears in \eqref{Shpot1} on the left. However, following the same logic as outlined above, by setting 
\begin{align} \label{Shpot2}
& m\rightarrow e-1-b \\
& n\rightarrow f-1-c \nonumber
\end{align}

we obtain a generalization of Shpot and Srivastava's result \cite[Eq. (31)] {Shpot&Sriv} (see \eqref{SSEq31}) in the form of a well known 3-part transformation, specifically \cite[Eq. 3.13.3(11)]{Luke}, and explicitly
\begin{align} \label{Shpot3}
& \displaystyle {\mbox{$_3$F$_2$}(a,b,c;\,e,f;\,1)}=\Gamma \left( 1-a \right) \Gamma \left( f \right) \Gamma \left( e \right)
\\  \nonumber
 &\times \left( {\frac {\Gamma \left( c-b \right) {\mbox{$_3$F$_2$}(b,-e+1+b,b+1-f;\,b-c+1,b+1-a;\,1)}}{\Gamma \left( c \right) \Gamma \left( e-b \right) \Gamma \left( b+1-a \right) \Gamma \left( -b+f \right) 
\mbox{}}} \right. \\ \nonumber
&\left. +{\frac {\Gamma \left( b-c \right) {\mbox{$_3$F$_2$}(c,-f+c+1,-e+c+1;\,c+1-a,c-b+1;\,1)}}{\Gamma \left( b \right) \Gamma \left( c+1-a \right) \Gamma \left( e-c \right) \Gamma \left( f-c \right) }} \right)\,. 
\end{align}

{\bf Remark 1.} Reversing \eqref{Shpot2}, gives a simple derivation of \cite[Eq. (31)] {Shpot&Sriv} from the known result \eqref{Shpot3}.\newline
{\bf Remark 2 and Digression.} In \eqref{Shpot1}, set $z=1$, sum (and reverse) the various series that arise and compare with \cite[Eq. (31)] {Shpot&Sriv} to find an interesting contiguity relation (also see \cite{Ebisu}):
\begin{equation}
\mapleinline{inert}{2d}{hypergeom([c, -n,b], [a, b+1+m], 1) = (-1)^m*GAMMA(1-b)*(sum((-1)^i*hypergeom([c, -n, b+i], [a, b+i+1], 1)/((-b-i)*GAMMA(i+1)*GAMMA(1-i+m)), i = 0 .. m))/GAMMA(-b-m)}{\[\displaystyle {\mbox{$_3$F$_2$}(c,-n,b;\,a,b+1+m;\,1)}={\frac { \left( -1 \right) ^{m}\Gamma \left( 1-b \right) }{\Gamma \left( -b-m \right) }\sum _{i=0}^{m}{\frac { \left( -1 \right) ^{i}{\mbox{$_3$F$_2$}(c,-n,b+i;\,a,b+i+1;\,1)}}{ \left( -b-i \right) \Gamma \left( i+1 \right) 
\mbox{}\Gamma \left( 1-i+m \right) }}}\]}\,.
\label{Shpot4}
\end{equation}  

\subsection{Obvious reduction from a $_4F_3$}
Set a pair of the top and bottom parameters in \eqref{Miller1} to be equal (e.g. $a=e$), and thereby trivially produce the three-part transformations among $_3F_2(1)$ reproduced above \eqref{Shpot3}.

With a similar choice of reduction ($e=a$) and reassignment of variable names in \eqref{GF_Imag}, another three-part transformation among $_3F_2(1)$ results. This result is equivalent to applying \eqref{Ginv} in the case $m=1,n=1,p=q=3$. Although this transformation involves terms clearly embedded in \cite[Tables 3.3, 3.4 and 3.5]{Luke} it is  given explicitly below:


\begin{align} \label{F32_Imag} 
&\displaystyle {\mbox{$_3$F$_2$}(a,b,c;\,f,e;\,1)}=-\frac {\Gamma \left( 1-b \right) \Gamma \left( 1-c \right) \Gamma \left( 1-a \right) }{\Gamma \left( 1-f \right) \Gamma \left( 1-e \right) }
\\ \nonumber
&\times \left( {\frac {\Gamma \left( -1+f \right) \Gamma \left( -e+f \right) {\mbox{$_3$F$_2$}(1+a-f,1+b-f,1+c-f;\,2-f,1+e-f;\,1)}}{\Gamma \left( f-a \right) \Gamma \left( -b+f \right) 
\mbox{}\Gamma \left( f-c \right) }}
\right. \\
&\left.+{\frac {\Gamma \left( -1+e \right) \Gamma \left( e-f \right) {\mbox{$_3$F$_2$}(1+a-e,1+b-e,1+c-e;\,2-e,1+f-e;\,1)}}{\Gamma \left( e-a \right) \Gamma \left( e-b \right) \Gamma \left( e-c \right) }}
\mbox{} \right)\,. \nonumber
\end{align}\newline

Because of the asymmetry among the parameters of the $_4F_3$ appearing in \eqref{Miller2}, various reductions from a $_4F_3$ to a $_3F_2$ can be made by choosing to equate different combinations of top and bottom parameters to yield different transformations. For example, the case $g=d$ produces a symmetric permutation equivalent to \eqref{Shpot3} whereas choosing $e=a$ and reassigning variables generates the following mixture, all of which are included (see Appendix C) in \cite[Tables (3.3), (3.4) and (3.5)]{Luke}:
\begin{align} \label{Miller2b} 
\displaystyle {\mbox{$_3$F$_2$}(a,b,c;\,e,f;\,1)}=& \mbox{}\Gamma \left( 1+b-e \right) \Gamma \left( 1-a \right) \Gamma \left( 1-c \right) \Gamma \left( f \right) 
\\ \nonumber
&\times\left( {\frac {{\mbox{$_3$F$_2$}(b,1+b-e,1+b-f;\,1+b-c,1+b-a;\,1)}}{\Gamma \left( 1-e \right) \Gamma \left( f-b \right) \Gamma \left( 1+b-c \right) \Gamma \left( 1+b-a \right) }}
\right. \\
&\left. +{\Gamma \left( e \right)\frac { {\mbox{$_3$F$_2$}(1+b-e,1+a-e,1+c-e;\,2-e,1-e+f;\,1)}}{\Gamma \left( 2-e \right) \Gamma \left( -c+e \right) \Gamma \left( e-a \right) \Gamma \left( 1-e+f \right) 
\mbox{}\Gamma \left( b \right) }} \right) \,. \nonumber
\end{align}

\subsection{Curious reduction from a $_4F_3(1)$}
A curious combination arises by setting $g=d$ in \eqref{Miller1} to obtain
\begin{align} \label{Mix1}
&\displaystyle {\mbox{$_3$F$_2$}(a,b,c;\,e,f;\,1)}=\frac {\Gamma \left( e \right) \Gamma \left( f \right) 
\mbox{}}{\sin \left( \pi\,d \right) } \\ \nonumber
& \times \left( -{\frac {\sin \left( \pi\, \left( a-d \right)  \right) {\mbox{$_3$F$_2$}(a,1+a-e,1+a-f;\,1+a-c,1+a-b;\,1)}\Gamma \left( c-a \right) 
\mbox{}\Gamma \left( b-a \right) }{\Gamma \left( f-a \right) \Gamma \left( e-a \right) \Gamma \left( c \right) \Gamma \left( b \right) }} \right. \\ \nonumber 
& \left.  -{\frac {\Gamma \left( a-b \right) \sin \left( \pi\, \left( b-d \right)  \right) {\mbox{$_3$F$_2$}(b,1+b-e,1+b-f;\,1+b-a,1+b-c;\,1)}\Gamma \left( c-b \right) }{\Gamma \left( a \right)\Gamma \left( c \right) 
\mbox{}\Gamma \left( f-b \right) \Gamma \left( e-b \right) }}
\right.  \\ \nonumber
&  \left. 
-{\frac {\Gamma \left( a-c \right) \Gamma \left( b-c \right) \sin \left( \pi\, \left( -d+c \right)  \right) {\mbox{$_3$F$_2$}(c,1+c-e,1+c-f;\,1+c-a,1+c-b;\,1)}}{\Gamma \left( a \right)\Gamma \left( f-c \right) \Gamma \left( e-c \right) 
\mbox{}\Gamma \left( b \right) }} \right)\,.   
\end{align}
Notice the appearance of the arbitrary parameter $d$ on the right-hand side but not on the left. Various possibilities were studied for various choices of the parameter $d$; even the choice $d=0$ yielded nothing other than one of the results quoted above. A typical choice of $d=1/2$ yields nothing more than a four-part transformation among $_3F_2(1)$ functions (but also see \eqref{P3}). However, a second even more curious transformation exists: set $g=a$ in \eqref{Miller2} to obtain another transformation, again with a superflous and arbitrary parameter $a$ appearing on the right-hand side, but not the left, this time embedded within the parameters of the hypergeometric functions themselves, rather than in a multiplicative factor. That is, after redefining parameters (including $a:=d$):

\begin{align} \label{Mix2}
&\displaystyle {\mbox{$_3$F$_2$}(a,b,c;\,e,f;\,1)}=\frac {\Gamma \left( 1+b-d \right) \Gamma \left( 1-c \right) \Gamma \left( 1-a \right) \Gamma \left( f \right) \Gamma \left( e \right) }{\Gamma \left( 1-d \right) } \\ \nonumber
& \times \left( {\frac {{\mbox{$_3$F$_2$}(b,1+b-e,1+b-f;\,1+b-a,1+b-c;\,1)}\Gamma \left( d-b \right) }{\Gamma \left( d \right) \Gamma \left( 1+b-a \right) 
\mbox{}\Gamma \left( 1+b-c \right) \Gamma \left( f-b \right) \Gamma \left( e-b \right) }} \right. \\ \nonumber
&\left. \hspace{.5cm}  -{\frac {\Gamma \left( d-b \right) {\mbox{$_4$F$_3$}(1,d,1-e+d,1-f+d;\,1+d-a,1+d-b,1+d-c;\,1)}}{\Gamma \left( f-d \right) \Gamma \left( e-d \right) \Gamma \left( 1+d-b \right) \Gamma \left( 1+d-c \right) 
\mbox{}\Gamma \left( 1+d-a \right) \Gamma \left( b \right) }} \right. \\  \nonumber
&\left. \hspace{.5cm}  -{\frac {\Gamma \left( d-1 \right) {\mbox{$_4$F$_3$}(1,1+a-d,1+b-d,1+c-d;\,2-d,1+e-d,1+f-d;\,1)}}{\Gamma \left( d \right) \Gamma \left( d-a \right) \Gamma \left( d-c \right) \Gamma \left( 1+e-d \right) 
\mbox{}\Gamma \left( 1+f-d \right) \Gamma \left( b \right) }}
\mbox{} \right) 
\end{align}

Although the case $d=0$ appears to be superficially interesting, after simplification, the result yields nothing more profound than the simple and well known reduction
\begin{equation}
\mapleinline{inert}{2d}{hypergeom([1, a, b, c], [2, e, f], 1) = (e-1)*(f-1)*(hypergeom([-1+a, -1+b, -1+c], [e-1, f-1], 1)-1)/((-1+a)*(-1+b)*(-1+c))}{\[\displaystyle {\mbox{$_4$F$_3$}(1,a,b,c;\,2,e,f;\,1)}={\frac { \left( e-1 \right)  \left( f-1 \right)  \left( {\mbox{$_3$F$_2$}(a-1,b-1,c-1;\,e-1,f-1;\,1)}-1 \right) }{ \left(a -1 \right) 
\mbox{} \left( b-1 \right)  \left(c -1 \right) }}\]}
\label{WellKnown}
\end{equation}

 The result \eqref{Mix2} however offers potential for investigating the properties of sums involving digamma functions, because the left hand-side vanishes under the action of the operator $\frac{\partial}{\partial\,d}$.

\section{Appendix: A primer on 3-part relations among $_3F_2(1)$} \label{sec:AppendixC}

The following should be read in conjunction with Section 3.13.3 of Luke's book \cite{Luke}; Slater's book \cite[Section 4.3.2]{Slater} also covers much of the same material.\newline

The three part relations among $_3F_2(1)$ were classified by Whipple (1923), summarized by Bailey \cite[Section 3.5]{Bailey1935} and reproduced by Luke \cite[Section 3.13.3]{Luke} and Slater \cite[Section 4.3.2]{Slater}. To say that the notation is obscure would be an understatement, and perhaps that is why these relations are often overlooked in the literature (e.g. see \eqref{Eq1p7b} and \eqref{Shpot3}). The following is an addendum to, and clarification of, Whipple's classification, the basis of which is a set of six parameters $r_i$ where $i=0(1)5$. To the best of my knowledge, these have never been explicitly listed in the literature. Working backwards from \cite[Table 3.3]{Luke} it is possible to obtain these parameters in terms of the top and bottom parameters of a canonical function $_3F_2(a,b,c;e,f;1)$ and thereby clarify the underlying algorithm. In summary,

\begin{align} \label{RRs} 
&\displaystyle {\it r_0} =5/6+c/3+b/3-2/3\,e-2/3\,f+a/3 \\ \nonumber
&\displaystyle {\it r_1}=-2/3\,b+e/3-1/6+f/3+a/3-2/3\,c \\ \nonumber
&\displaystyle {\it r_2}=f/3-2/3\,a+e/3-1/6+b/3-2/3\,c \\ \nonumber
&\displaystyle {\it r_3}=c/3-1/6+e/3-2/3\,b+f/3-2/3\,a \\ \nonumber
&\displaystyle {\it r_4}=e/3-1/6+b/3-2/3\,f+a/3+c/3 \\ \nonumber
&\displaystyle {\it r_5}=-1/6+b/3-2/3\,e+f/3+a/3+c/3\,, \nonumber
\end{align} 
from which, with the help of \cite[Eq. (3.13.3(13)]{Luke} it is possible to reproduce \cite[Table 3.3]{Luke} in its entirety, and calculate representative (i.e. mixed combination of $a,b,c,e,f$ ) top and bottom parameters $\alpha_{lmn}$ and $\beta_{mn}$ respectively (see \cite[Eq. 3.13(13)]{Luke}), {\bf labelled} by distinct permutations of integers $l,m,n=0(1)5$. Whipple then introduces two fundamental functions $F_p(u;v,w)$ and $F_n(u;v,w)$ \cite[Eqs. 3.13.3(14) and 3.13.3(15)]{Luke} defined in terms of $_3F_2(a,b,c;e,f;1)$, where each of the labels $u,v,w$ take on one of the distinct, but different, numbers $0,1...5$. Representative independent mixtures of the top and bottom parameters are classified and listed in \cite[Tables 3.4 and 3.5]{Luke}; any missing combinations from the tables represent simple, irrelevant permutations among the top three, or between the bottom two, parameters of the canonical $_3F_2(1)$. The important point to note is that only the first parameter of $F_p(u;v,w)$ or $F_n(u;v,w)$ is important; the second and third parameters represent different mixtures among the top and bottom parameters and are conveniently omitted unless needed for clarity or specificity. Because of the two-part (Thomae) relations, all twenty of these functions with the same first parameter are equal; that is 

\begin{equation}
F_p(u;v_1,w_1)=F_p(u;v_2,w_2)
\end{equation}
and
\begin{equation}
F_n(u;v_1,w_1)=F_n(u;v_2,w_2)
\end{equation}
for any permitted combination of $u,v,w =0(1)5$ (no duplication). In the following, it is assumed that the arguments of any function on either side of an equality sign are such that the series representation converges (if $s\equiv e+f-a-b-c$, then $\Re(s)>0$); otherwise, the equality between sides must be interpreted in the sense of analytic continuation. With this notation, \cite[Eq. 3.13.3(11)]{Luke} (also \cite[Eq. 7.4.4(3)]{prudnikov} or \eqref{Shpot3}), a well known three-part relation among three particular $_3F_2(1)$, can be written in labelled form 

\begin{equation}
\displaystyle {\it F_p} \left( 0;4,5 \right) ={\frac {\pi\,\Gamma \left( \alpha _{ 023}  \right) }{\sin \left( \pi\,\beta _{23}  \right) }} \left( {\frac {{\it F_n} (2) }{\displaystyle\Gamma \left( \alpha_{134}  \right) \Gamma \left( \alpha _{135}\right) 
\mbox{}\Gamma \left( \alpha _{345} \right) }}-{\frac {{\it F_n} \left( 3 \right) }{\Gamma \left( \alpha_{124}  \right) \Gamma \left( \alpha _{125}  \right) \Gamma \left( \alpha_{245}\right) }}\,,
\mbox{} \right)
\label{Q3}
\end{equation}
or, in expanded form
\begin{align} \label{P3}
\displaystyle {\it F_p} \left( 0;4,5 \right)& =-\Gamma \left( c+1-b \right) \Gamma \left( b-c \right) \Gamma \left( 1-a \right)
\\ \nonumber
&\times \left( {\frac {{\it F_n} \left( 2;3,1 \right) }{\Gamma \left( -b+e \right) \Gamma \left( -b+f \right) \Gamma \left( c \right) }}-{\frac {{\it F_n} \left( 3;1,2 \right) }{\Gamma \left( e-c \right) \Gamma \left( f-c \right) \Gamma \left( b \right) }}
\mbox{} \right)\,. 
\end{align}
In \eqref{P3}, an arbitrary choice of (superfluous) second and third parameters of $F_n(2)$ and $F_n(3)$ have been included, and the arguments of the $\Gamma$ functions have been written explicitly in terms of the underlying parameters, to specify one of the 20 possibilities equivalent to \eqref{Shpot3}. By limiting the left-hand side to one particular (the canonical) $_3F_2(1)$, and removing the second and third parameters from the right-hand side, \eqref{Q3} represents a family of forty 3-part relations among different combinations of $_3F_2(1)$ selected by different combinations of second and third parameters (see \cite[Table 3.13.3.5]{Luke}) on the right-hand side. Three other (equivalent families of) three-part relations among $_3F_2(1)$ are known and listed in \cite[Eqs. 7.4.4(4) -(6)]{prudnikov}. The first two of these, once parsed according to the tables and rules cited above, and illustrated by one specific instance each, are

\begin{align} \label{P4}
& \displaystyle {\it F_p} \left( 0;4,5 \right) =  \Gamma \left( 1+c-e \right) \\ \nonumber 
&\times\,\left( {\frac {\Gamma \left( e-a-b \right) {\it F_n} \left( 5;0,3 \right) \Gamma \left( 1+a+b-e \right) }{\Gamma \left( s \right) \Gamma \left( -a+e \right) 
\mbox{}\Gamma \left( -b+e \right) }} 
 +{\frac {\Gamma \left( a+b-e \right) {\it F_n} \left( 3;0,5 \right) \Gamma \left( 1+e-b-a \right) }{\Gamma \left( b \right) \Gamma \left( f-c \right) \Gamma \left( a \right) }} \right)
\end{align} 
and
\begin{align} \label{P5}
\displaystyle {\it F_p} \left( 0;4,5 \right)  = \Gamma \left( 1-f+b \right) 
\mbox{}\Gamma \left( 1+c-f \right) \Gamma \left( 1-f+a \right)\, \\ \nonumber
\times \left( {\frac {{\it F_n} \left( 4;0,2 \right) }{\Gamma \left( s \right) \Gamma \left( f \right) \Gamma \left( 1-f \right) }}+{\frac {{\it F_p} \left( 5; 0,4 \right) }{\Gamma \left( b \right) \Gamma \left( c \right) \Gamma \left( a \right) }} 
\mbox{} \right)\,. 
\end{align}
Surprisingly, the third equation  \cite[Eq. 7.4.4(6)]{prudnikov} cited above, once parsed in the same manner turns out to be equivalent to \eqref{P5} with the change ${\it F_p} \left( 4;0,2 \right)= {\it F_p} \left( 4;2,3 \right) $ and is therefore not independent, although when expanded in its full glory, this is not evident.\newline

Finally, it is noted that three more independent equations between the $F_n$ and $F_p$ functions can be found by changing the signs of the $r_i$ terms in \cite[Eq. (13)]{Luke}. This has the effect of redefining the parameters  $\alpha_{lmn}\rightarrow 1-\alpha_{lmn}$ and $\beta_{mn}\rightarrow 2-\beta_{mn}$ in the various tables, as well as converting $F_p(u) \rightleftarrows F_n(u)$ and $s\rightarrow 1-s$. Thus \eqref{P3} becomes

\begin{align} \label{P3A}
\displaystyle {\it F_n} \left( 0 \right)& =-\frac {\pi\,\Gamma \left( a \right) }{\sin \left( \pi\, \left( 1+b-c \right)  \right) } \\ \nonumber
 & \times\,\left( {\frac {{\it F_p} \left( 2 \right) }{\Gamma \left( 1-e+b \right) \Gamma \left( 1-f+b \right) \Gamma \left( 1-c \right) }}
\mbox{}  
 -{\frac {{\it F_p} \left( 3 \right) }{\Gamma \left( 1+c-e \right) \Gamma \left( 1+c-f \right) \Gamma \left( 1-b \right) }} \right), 
\end{align}
this time omitting the superfluous second and third parameters from $F_p(0), F_p(2)$ and $F_p(3)$, yet retaining all combinations of parameters in a form that a patient reader could identify as $\alpha_{lmn}$ or $\beta_{mn}$ from Luke's Tables. Similarly, \eqref{P4} and \eqref{P5} become

\begin{align} \label{P4A}
\displaystyle {\it F_n} \left( 0 \right) &=\frac {\pi\,\Gamma \left( e-c \right) }{\sin \left( \pi\, \left( 1+e-b-a \right) 
\mbox{} \right) } \\ \nonumber
& \times \left( {\frac {{\it F_p} \left( 5 \right) }{\Gamma \left( 1-e+b \right) \Gamma \left( 1-e+a \right) \Gamma \left( 1-s \right) }}
\mbox{}-{\frac {{\it F_p} \left( 3 \right) }{\Gamma \left( 1-b \right) \Gamma \left( 1+c-f \right) \Gamma \left( 1-a \right) }} \right) 
\end{align}
and
\begin{align} \label{P5A}
\displaystyle {\it F_n} \left( 0 \right) = \Gamma \left( -b+f \right) \Gamma \left( f-c \right) \Gamma \left( f-a \right)\left( -{\frac {\sin \left( \pi\,f \right) {\it F_p} \left( 4 \right) }{\pi\,\Gamma \left( 1-s \right) }}+{\frac {{\it F_n} \left( 5 \right) }{\Gamma \left( 1-b \right) \Gamma \left( 1-c \right) \Gamma \left( 1-a \right) }}
\mbox{} \right) \,.
\end{align}

A collection consisting of six such equations is sufficient to interrelate all possible three-term relations among 120 $_3F_2(1)$ that can be identified as $F_p(u)$ and/or $F_n(u)$ \cite[Section 4.3.2]{Slater}. Luke \cite[Eq. 3.13(26)]{Luke} and Slater \cite[Eq. (4.3.2.5)]{Slater} go on to reproduce an example from Bailey that uses the above to relate $F_n(0), F_p(0)$ and $F_p(5)$ algebraically. Luke does not say which six equations he used. Similarly, at this same point, Slater refers to ``the relation corresponding to \cite[Eq. 4.3.2.1)]{Slater} which connects $F_p(5),F_n(0)$ and $F_n(2)$..." but she  never identifies the corresponding relation. Both references in Luke and Slater correspond to \eqref{P3} here. It is reasonable, but not assured, to assume that Slater, Luke (and Bailey) based this example on the three well-established 3-part relations \cite[Eqs. 7.4.4(3) -(5)]{prudnikov} that gave rise to the above (Note: Bailey refers to Hardy and Whipple in a footnote at this point but also does not say which equations were used); Slater \cite[Eq. (4.3.2.4)]{Slater} identifies one. How are other relations found?\newline

Dealing with the case of $F_n(1)$ and $F_p(1)$ (missing from all six of the above), note that $F_n(1;u,v)$ and $F_n(2;u,v)$ are related by an interchange of two of $a,b$ or $c$. For example, under the (symmetric) exchange $a\leftrightarrows b$, \eqref{P3} becomes

 
\begin{align} \label{P3C}
\displaystyle {\it F_p} \left( 0 \right) ={\frac {\pi\,\Gamma \left( 1-b \right) }{\sin \left( \pi\, \left( a-c \right)  \right) } \left( -{\frac {{\it F_n} \left( 1 \right) }{\Gamma \left( -a+e \right) \Gamma \left( f-a \right) \Gamma \left( c \right) }}
\mbox{}+{\frac {{\it F_n} \left( 3 \right) }{\Gamma \left( e-c \right) \Gamma \left( f-c \right) \Gamma \left( a \right) }} \right) }
\end{align}

along with its complement
\begin{align} \label{P3Cc}
&\displaystyle {\it F_n} \left( 0 \right)=\\ \nonumber
&\hspace{0.5cm}\frac {\pi\,\Gamma \left( b \right) }{\sin \left( \pi\, \left( a-c \right)  \right) } 
\left( {\frac {{\it F_p} \left( 1 \right) }{\Gamma \left( 1-e+a \right) \Gamma \left( 1-f+a \right) \Gamma \left( 1-c \right) }}
\mbox{}-{\frac {{\it F_p} \left( 3 \right) }{\Gamma \left( 1+c-e \right) \Gamma \left( 1+c-f \right) \Gamma \left( 1-a \right) }} \right)  \,.
\end{align}
and this establishes that all the basic functions $F_n$ and $F_p$ are at least accessible from the independent relations. When written in labelled form \eqref{P3C} becomes
\begin{equation}
\displaystyle {\it F_p} \left( 0 \right) ={\frac {\pi\,\Gamma \left( {\it \alpha} _{ 013} \right)  }{\sin
 \left( \pi\,\beta_{13}  \right) }} \left( {\frac {{\it F_n} \left( 1 \right) }{\Gamma \left( {\it \alpha} _{234}  \right) \Gamma \left( {\it \alpha} _{235}  \right) \Gamma \left( {\it \alpha} _{345}  \right) 
\mbox{}}}-{\frac {{\it F_n} \left( 3 \right) }{\Gamma \left( {\it \alpha} _{214}  \right) \Gamma \left( {\it \alpha} _{215}   \right) \Gamma \left( {\it \alpha}_{145}  \right) }}\,.
\mbox{} \right)
\label{P3Ca}
\end{equation}
Notice  that \eqref{Q3} and \eqref{P3Ca}, related by the interchange $a \leftrightarrows b$, are equivalently related by the interchange of the numeric labels $1\leftrightarrows2$; this demonstrates the fundamental basis of the notation - {\it the equality of the six equations chosen as a basis is invariant under interchange of numeric labels.} With this understanding, all other relationships can be found. For example, consider the independently derived result \eqref{Miller2b} 
\begin{equation}
\mapleinline{inert}{2d}{fp(0) = (fn(2)*sin(Pi*e)/(Pi*GAMMA(-b+f))+fp(4)/(GAMMA(e-c)*GAMMA(-a+e)*GAMMA(b)))*GAMMA(1-c)*GAMMA(1-a)*GAMMA(1-e+b)}{\[\displaystyle {\it F_p} \left( 0 \right) = \left( {\frac {\sin \left( \pi\,e \right){\it F_n} \left( 2 \right)  }{\pi\,\Gamma \left( f-b \right) }}+{\frac {{\it F_p} \left( 4 \right) }{\Gamma \left( e-c \right) \Gamma \left( -a+e \right) \Gamma \left( b \right) }}
\mbox{} \right) \Gamma \left( 1-c \right) 
\mbox{}\Gamma \left( 1-a \right) \Gamma \left( 1-e+b \right) \]}
\label{B8}
\end{equation}
along with its complement
\begin{equation}
\mapleinline{inert}{2d}{fn(0) = (fn(4)/(GAMMA(1+c-e)*GAMMA(1-e+a)*GAMMA(1-b))-fp(2)*sin(Pi*e)/(GAMMA(1-f+b)*Pi))*GAMMA(c)*GAMMA(a)*GAMMA(-b+e)}{\[\displaystyle {\it F_n} \left( 0 \right) = \left( {\frac {{\it F_n} \left( 4 \right) }{\Gamma \left( 1+c-e \right) \Gamma \left( 1+a-e \right) \Gamma \left( 1-b \right) }}-{\frac {\sin \left( \pi\,e \right){\it F_p} \left( 2 \right)  }{\pi\,\Gamma \left( 1-f+b \right) }}
\mbox{} \right) 
\mbox{}\Gamma \left( c \right) \Gamma \left( a \right) \Gamma \left( -b+e \right) \]}\,.
\label{B8c}
\end{equation}

To demonstrate that \eqref{B8} can be obtained from one (or more) of the six independent relations given above, start with \eqref{P5}, written in labelled form

\begin{equation}
\displaystyle {\it F_p} \left( 0 \right) =\left({\frac {\sin \left( \pi\,\beta _{50}  \right)
 {\it F_n} \left( 4 \right) }{\Gamma \left( \alpha _{123}   \right) 
\mbox{}\pi}}+{\frac {  {\it F_p} \left( 5 \right) }{\Gamma \left( \alpha _{245}\right) 
\Gamma \left( \alpha_{345} \right) \Gamma \left( \alpha _{145}  \right) }}\right) \Gamma \left( \alpha_{024} \right) \Gamma \left( \alpha_{ 014}  \right) \Gamma \left( \alpha _{034} \right)\,.
\label{P5a}
\end{equation}
Interchange all labels $4\leftrightarrows2$ followed by the interchange $5\leftrightarrows4$ to find
\begin{equation}
\displaystyle {\it F_p} \left( 0 \right) = \left( {\frac {\sin \left( \pi\,\beta _{40}  \right) {\it F_n} \left( 2 \right) }{\Gamma \left( \alpha _{153}  \right) 
\mbox{}\pi}}+{\frac {{\it F_p} \left( 4 \right) }{\Gamma \left( \alpha _{524} \right) \Gamma \left( \alpha _{324}\right) \Gamma \left( \alpha _{124}  \right) }} \right) 
\mbox{}\Gamma \left( \alpha _{012}  \right) \Gamma \left( \alpha_{052}  \right) \Gamma \left( \alpha _{ 032}  \right) 
\label{B8a}
\end{equation}
which, when written in terms of the basic parameters, is identifiable as \eqref{B8}.\newline

To derive \eqref{F32_Imag} is more interesting. Written in this notation,  \eqref{F32_Imag} becomes
\begin{align} \label{B7}
\displaystyle {\it F_p} \left( 0 \right)& =\frac {\sin \left( \pi\,e \right) \sin \left( \pi\,f \right) \Gamma \left( 1-a \right) \Gamma \left( 1-b \right)
\mbox{}\Gamma \left( 1-c \right) }{\sin \left( \pi\, \left( e-f \right)  \right) } \\ \nonumber
&\times  \left( {\frac {{\it F_p} \left( 4 \right) }{\sin \left( \pi\,e \right) \Gamma \left( -b+e \right) \Gamma \left( e-c \right) \Gamma \left( -a+e \right) }}
\mbox{}-{\frac {{\it F_p} \left( 5 \right) }{\sin \left( \pi\,f \right) \Gamma \left( -b+f \right) \Gamma \left( f-c \right) \Gamma \left( f-a \right) }} \right) 
\end{align}
and its derivation from any one of the six independent relations by interchanging labels becomes questionable, because none of them individually relate three $F_p$ functions. Choose \eqref{P3A} and \eqref{P4A} as a convenient staring point, equate the right-hand sides of both and solve for $F_p(2)$. Written in labelled notation, the solution is
\begin{align} \nonumber
\displaystyle {\it F_p} \left( 2 \right) =&\left( \frac{1}{ \Gamma \left( \alpha_{035}  \right)  }+{\frac {\Gamma \left( \alpha _{124}\right) \sin \left( \pi\,\beta_{23}  \right) }{\Gamma \left( \alpha _{145}  \right) 
\Gamma \left( \alpha _{023}  \right) \sin \left( \pi\,\beta_{35}  \right) }}
\mbox{} \right){\frac {  \Gamma \left( \alpha _{ 012}  \right) \Gamma \left( \alpha _{024}  \right) 
\mbox{}\Gamma \left( \alpha _{ 025}  \right)  }{\Gamma \left( \alpha _{034}  \right) \Gamma \left( \alpha _{013}  \right) }}F_{{p}} \left( 3 \right)
\mbox{}\\&
-F_{{p}} \left( 5 \right){\frac { \Gamma \left( \alpha _{124}  \right) \sin \left( \pi\,\beta_{23}  \right) \Gamma \left( \alpha _{024}  \right)
\Gamma \left( \alpha _{ 012}  \right) }{\Gamma \left( \alpha _{ 045}  \right) \sin \left( \pi\,\beta _{35}  \right) \Gamma \left( \alpha _{ 015}  \right)
\Gamma \left( \alpha _{145}  \right) }}\,.
\label{EqY}
\end{align}
\newline
Now perform the interchange of numeric labels $3\leftrightarrows4$, followed by the interchange $0\leftrightarrows2$, revert to the representation in terms of the underlying parameters, and after some simplification, \eqref{B7} will be found.\newline

In the section dealing with the generalization of Shpot and Srivastava's result \eqref{Shpot1}, it was claimed that their \cite[Eq.(31)]{Shpot&Sriv} (see \eqref{SSEq31}) could easily be obtained from the ``well known" result \eqref{Shpot3}. A quick scan of \cite[Table 3.5]{Luke} shows that the two $_3F_2(1)$ appearing in \eqref{Shpot3} can be identified as $F_n(2;3,1)$ and $F_n(3;1,2)$ as discussed above. Simple substitution into \eqref{P3} will yield \eqref{Shpot3}, justifying the remark that it is at least a ``known", if not a ``well known", result. In the case of Paris and Miller, their extended result \eqref{Eq1p7b} once parsed as discussed here can be written

\begin{equation}
\mapleinline{inert}{2d}{fp(0) = (sin(Pi*e)*fn(3)/(Pi*GAMMA(f-c))+fp(4)/(GAMMA(-b+e)*GAMMA(-a+e)*GAMMA(c)))*GAMMA(1+c-e)*GAMMA(1-b)*GAMMA(1-a)}{\[\displaystyle {\it F_p} \left( 0 \right) = \left( {\frac {\sin \left( \pi\,e \right) {\it F_n} \left( 3 \right) }{\pi\,\Gamma \left( f-c \right) }}+{\frac {{\it F_p} \left( 4 \right) }{\Gamma \left( -b+e \right) \Gamma \left( -a+e \right) \Gamma \left( c \right) }}
\mbox{} \right) \Gamma \left( 1+c-e \right) 
\mbox{}\Gamma \left( 1-b \right) \Gamma \left( 1-a \right) \]}
\label{P_and_M}
\end{equation}

where specifically $F_p(0)=F_p(0;4,5), F_p(4)=F_p(4;1,2)$ and $F_n(3)=F_n(3;1,2)$. To derive \eqref{P_and_M} from the above, solve for $F_n(2)$ in \eqref{P3} and substitute into \eqref{B8}. So in the sense discussed here \eqref{P_and_M} is a known result.

\section{A collection of $_3F_2(1)$, $_4F_3(1)$ and some lemmas}\label{sec:Misc}

The following is a collection of relevant results gathered from sources scattered throughout the literature, plus a few lemmas.\newline

\begin{itemize}
\item{\bf Minton Karlsson } \label{MinK}

Minton \cite{Minton} and Karlsson \cite{Karlsson} show that, when a top parameter exceeds a bottom parameter by a positive integer $n$, the order of any $_{p+1}F_p(1)$ can be reduced by one and replaced by a sum of $n$ terms. In the case $p=3$, this gives 
\begin{align} \label{Mk}
&\displaystyle \mbox{$_4$F$_3$}(a,b,c,e+n;\,e,f,g;\,1)= \\ \nonumber
&\hspace{.35cm}{\frac {\Gamma \left( 1+n \right) \Gamma \left( e \right) \Gamma \left( f \right) \Gamma \left( g \right) }{\Gamma \left( a \right) \Gamma \left( b \right) 
\mbox{}\Gamma \left( c \right) }\sum _{k=0}^{n}{\frac {\Gamma \left( k+a \right) \Gamma \left( b+k \right) \Gamma \left( c+k \right) {\mbox{$_3$F$_2$}(k+a,b+k,c+k;\,f+k,g+k;\,1)}}{\Gamma \left( k+1 \right) 
\mbox{}\Gamma \left( 1+n-k \right) \Gamma \left( e+k \right) \Gamma \left( f+k \right) \Gamma \left( g+k \right) }}}.
\end{align}\newline
See also \cite{Gottschalk}.

\item{\bf Sheppard-Andersen Theorem} \label{ShepAnd}

Based on \cite[Corollary 3.3.4]{AndAskRoy}, the general result for a k-balanced, terminating $_3F_2(1)$ is given by
\begin{align} \label{HRKX}
&\displaystyle {\mbox{$_3$F$_2$}(a,b,-n;\,c,k-n+a+b-c;\,1)}=\\ \nonumber
&-\frac {\pi\, \left( -1 \right) ^{n}\Gamma \left( a-c+1 \right)  \left( -1 \right) ^{k}
\mbox{}\Gamma \left( k-n+a+b-c \right) \Gamma \left( n+1 \right) \Gamma \left( k \right) \Gamma \left( c \right) }{\sin \left( \pi\, \left( c-b \right)  \right) \Gamma \left( k-n+b-c \right) 
\mbox{}\Gamma \left( k+b-c+a \right) \Gamma \left( n+c \right) \Gamma \left( a \right) }\\ \nonumber
&\hspace{2cm}\times\sum _{j=0}^{N}{\frac {\Gamma \left( a+j \right) }{\Gamma \left( 1+n-j \right) \Gamma \left( -b+c-k+1+j \right) 
\mbox{}\Gamma \left( a-n-c+1+j \right) \Gamma \left( 1+j \right) \Gamma \left( k-j \right) }}
\end{align}

where $N=min(k-1,n)$. In the case $k=2$ a simpler result \cite[Eq. (2.9), misprinted]{RoyAmMath} is

\begin{align} \label{SARoy}
&\displaystyle {\mbox{$_3$F$_2$}(a,b,-n;\,c,2-n+a+b-c;\,1)} \\ \nonumber
&\hspace{.9cm}={\frac {\Gamma \left( -b+c+n-1 \right) 
\mbox{}\Gamma \left( n+c-a \right) \Gamma \left( c-a-b-1 \right) \Gamma \left( c \right) }{\Gamma \left( c-b-1 \right) \Gamma \left( c-a \right) \Gamma \left( -b+c+n-1-a \right) 
\mbox{}\Gamma \left( n+c \right) } \left( 1-{\frac {na}{ \left( c-b-1 \right)  \left( -a+n+c-1 \right) }} \right) }
\end{align}

{\bf Remark:} The result \eqref{HRKX} is usually referenced in the literature to an inaccessible paper by Sheppard \cite{Sheppard}, where it is also usually noted that Andersen \cite{Andersen} obtained \eqref{SARoy} independently. In  fact, Andersen obtained the following result (transcribed in hypergeometric notation)
\begin{align} \label{And1}
\displaystyle \mbox{$_3$F$_2$}&(1,a,m+b;\,2+m,a+b;\,1)\\ \nonumber
&={\frac { \left( -1 \right) ^{m}\sin \left( \pi\,a \right) \sin \left( \pi\,b \right) \Gamma \left( 1-b-m \right) \Gamma \left( b+a \right) \Gamma \left( -a+1 \right) 
\mbox{}\Gamma \left( 2+m \right) }{ \left( -m+a-1 \right) {\pi}^{2} \left( b-1 \right) }}-{\frac { \left( 1+m \right)  \left( b+a-1 \right) }{ \left( -m+a-1 \right)  \left( b-1 \right) }}
\end{align}
and a similar result for the case where the bottom parameter is $a+b+1$, both of which are special cases of \cite[Entry 28]{Milgram447}.


\item{\bf Whipples tranformation} \label{Whipple1}

The result \eqref{HRKX} is based upon the following transformation of a terminating 1-balanced $_4F_3(1)$ due to Whipple (1926) (see \cite[Theorem 3.3.3]{AndAskRoy})

\begin{align} \label{WhipXf1}
&\displaystyle {\mbox{$_4$F$_3$}(a,b,c,-n;\,d,e,f;\,1)}\\ \nonumber
&={\frac {\Gamma \left( n+e-a \right) \Gamma \left( n+f-a \right) \Gamma \left( e \right) \Gamma \left( f \right) {\mbox{$_4$F$_3$}(a,-n,d-b,d-c;\,d,a+1-n-f,a+1-n-e;\,1)}
\mbox{}}{\Gamma \left( e-a \right) \Gamma \left( f-a \right) \Gamma \left( n+e \right) \Gamma \left( n+f \right) }}
\end{align}
where $f=a+b+c+1-d-e-n$.\newline

Suppose $d=c-1$. Then \eqref{WhipXf1} gives the following result for a special (Minton-Karlsson), 1-balanced, terminating $_4F_3(1)$

\begin{align} \label{WhipXf2}
\displaystyle \mbox{$_4$F$_3$} &\left( -n,a,b,c;c-1,e,a+b+2-n-e;1 \right)\\ \nonumber
& =A{\frac {\Gamma \left( a+b+2-n-e \right) 
\mbox{}\Gamma \left( n+e-a-1 \right) \Gamma \left( b+1-e \right) \Gamma \left( e \right) }{\Gamma \left( b+2-n-e \right) \Gamma \left( a+b+2-e \right) \Gamma \left( n+e \right) 
\mbox{}\Gamma \left( e-a \right) }}
\end{align}
where 
\begin{equation}
\displaystyle A=-{e}^{2}+ \left( a+b-n+2 \right) e-{\frac { \left( b+1 \right)  \left( a-c+1 \right)  \left( n-1 \right) }{c-1}}+{\frac { \left( -b-1+n \right) ac}{c-1}} \,.
\label{WhipA}
\end{equation}
Other variations are apparent.

\item{\bf A special {\mbox{$_{q+1}$F$_q$}}}
 
Prudnikov et. al. \cite[Eq. 7.10.2(6)]{prudnikov} give the following general result
\begin{equation} \label{pFq}
\displaystyle {\mbox{$_{q+1}$F$_q$}(a,b,\dots,b;\,b+1,\dots,b+1;\,1)}={\frac { \left( -1 \right) ^{q-1}{b}^{q}\Gamma \left( 1-a \right) 
\mbox{}}{ \left( q-1 \right) !}{\frac {\partial ^{q-1}}{\partial {b}^{q-1}}} \left( {\frac {\Gamma \left( b \right) }{\Gamma \left( 1+b-a \right) }} \right) }.
\end{equation}

When $q=4$ and $b=1$, we find
\begin{align} \label{5F4b1}
\displaystyle  \left( 1-a \right) &\mbox{$_5$F$_4$}(1,1,1,1,a;\,2,2,2,2;\,1)=1/6\, \left( \Psi \left( 2-a \right) +\gamma \right) ^{3}\\ \nonumber
&+ \left( 1/12\,{\pi}^{2}-1/2\,\Psi^{\prime} \left(2-a \right)  \right)  \left( \Psi \left( 2-a \right) +\gamma \right) 
\mbox{} 
+1/3\,\zeta \left( 3 \right) +1/6\,\Psi^{(2)} \left(2-a \right) 
\end{align}

and similarly with $q=5$,
\begin{align} \label{6F5b1}
&\displaystyle  \left( 1-a \right) {\mbox{$_6$F$_5$}(1,1,1,1,1,a;\,2,2,2,2,2;\,1)}=1/8\, \Psi^{\prime} \left( 2-a \right)^{2}+{\frac {{\pi}^{4}}{160}}-1/24\,\Psi^{(3)} \left(2-a \right)  
\\ \nonumber 
&+1/24\, \left( \Psi \left( 2-a \right) +\gamma \right) ^{4}
\mbox{}+ \left( 1/3\,\zeta \left( 3 \right) 
+1/6\,\Psi^{(2)} \left(2-a \right)  \right)  \left( \Psi \left( 2-a \right) +\gamma \right) \\ \nonumber
& - \left( 1/4\, \left( \Psi \left( 2-a \right) +\gamma \right) ^{2}+1/24\,{\pi}^{2} \right) \Psi^{\prime} \left(2-a \right)+1/24\,{\pi}^{2} \left( \Psi \left( 2-a \right) +\gamma \right) ^{2}.
\end{align}

\item{\bf A result from Shpot and Srivastava \cite{Shpot&Sriv}}

The following reproduces \cite[Eq. (31)]{Shpot&Sriv} 

\begin{equation}
\mapleinline{inert}{2d}{Eq31 := hypergeom([a, b, c], [b+1+m, c+1+n], 1)*GAMMA(b)*GAMMA(c)/(GAMMA(b+1+m)*GAMMA(c+1+n)) = T(a, b, c, m, n)+T(a, c, b, n, m)}{\[\displaystyle {\frac {{\mbox{$_3$F$_2$}(a,b,c;\,b+1+m,c+1+n;\,1)}\Gamma \left( b \right) \Gamma \left( c \right) }{\Gamma \left( b+1+m \right) \Gamma \left( c+1+n \right) }}=T \left( a,b,c,m,n \right) 
\mbox{}+T \left( a,c,b,n,m \right) \]}
\label{SSEq31}
\end{equation}
where
\begin{equation}
\mapleinline{inert}{2d}{T(a, b, c, m, n) = GAMMA(b)*GAMMA(1-a)*GAMMA(c-b)*hypergeom([b, -m, b-c-n], [1+b-a, 1+b-c], 1)/(GAMMA(1+b-a)*GAMMA(n+1+c-b)*GAMMA(m+1))}{\[\displaystyle T \left( a,b,c,m,n \right) ={\frac {\Gamma \left( b \right) \Gamma \left( 1-a \right) \Gamma \left( c-b \right) {\mbox{$_3$F$_2$}(b,-m,b-c-n;\,1+b-a,1+b-c;\,1)}}{\Gamma \left( 1+b-a \right)
\mbox{}\Gamma \left( n+1+c-b \right) \Gamma \left( m+1 \right) }}\]}
\label{Tabcmn}
\end{equation}

\item{\bf From a previous work \cite[Theorem 2.2]{MilgDigSums}}
\begin{align} \label{Eq2p2}
\displaystyle {\mbox{$_4$F$_3$}(1,1,1,1;\,2,2,d;\,1)}=& \left( d-1 \right)  \bigg{(}\Psi^{(2)} \left(d-1 \right)/2 - \left( \Psi^{\prime} \left(d-1 \right) +{\pi}^{2}/6\, \right)  \left( \gamma+\Psi \left( d-1 \right)  \right) 
\mbox{} \\ \nonumber
&\left. +2\,\zeta \left( 3 \right)  -{\frac {1}{\Gamma \left( 2-d \right) }\sum _{k=0}^{\infty }{\frac {\Psi^{\prime} \left(k+2 \right) \Gamma \left( 3-d+k \right) }{\Gamma \left( 1+k \right)  \left( 1+k \right) ^{2}}}}
\mbox{} \right)
\end{align}

\item{\bf Lemma 1}

Because
\begin{equation}
\mapleinline{inert}{2d}{Sum(GAMMA(n+k-b)*(-x)^k/(GAMMA(k+1)^2*GAMMA(n-k)), k = 0 .. n-1) = [GAMMA(n-b)*hypergeom([n-b, -n+1], [1], x)/GAMMA(n) = -(-1)^n*GAMMA(n-b)^2/(GAMMA(-b+1)*GAMMA(n)^2)]}{\[\displaystyle \sum _{k=0}^{n-1}{\frac {\Gamma \left( n+k-b \right)  \left( -1 \right) ^{k}}{  \Gamma \left( k+1 \right)   ^{2}\Gamma \left( n-k \right) }}={\frac {\Gamma \left( n-b \right) {\mbox{$_2$F$_1$}(n-b,-n+1;\,1;\,1)}
\mbox{}}{\Gamma \left( n \right) }}=-{\frac { \left( -1 \right) ^{n} \Gamma \left( n-b \right)  ^{2}}{\Gamma \left( -b+1 \right) \Gamma \left( n \right)  ^{2}}}\]}\,,
\label{Cxs1}
\end{equation}

differentiate with respect to $b$ giving
\begin{equation}
\mapleinline{inert}{2d}{Sum(GAMMA(n+k-b)*(-1)^k*Psi(n+k-b)/(GAMMA(k+1)^2*GAMMA(n-k)), k = 0 .. n-1) = -(-1)^n*Pi*(-Psi(1-b)+2*Psi(n-b))*GAMMA(b)/(sin(Pi*b)*GAMMA(n)^2*GAMMA(-n+b+1)^2)}{\[\displaystyle \sum _{k=0}^{n-1}{\frac {\Gamma \left( n+k-b \right)  \left( -1 \right) ^{k}\Psi \left( n+k-b \right) }{ \Gamma \left( k+1 \right)   ^{2}
\mbox{}\Gamma \left( n-k \right) }}=-{\frac { \left( -1 \right) ^{n}\pi\, \left( -\Psi \left( 1-b \right) +2\,\Psi \left( n-b \right)  \right) \Gamma \left( b \right) }{\sin \left( \pi\,b \right) 
\mbox{}  \Gamma \left( n \right) ^{2}  \Gamma \left( -n+b+1 \right)  ^{2}}}\]}\,,
\label{CxSum}
\end{equation}

\item{\bf Lemma 2}

From \cite[Entry 13]{Milgram447}, corresponding to ${\mbox{$_3$F$_2$}(1,1,a;\,1+a,1+a;\,1)}$, a case contiguous to Whipple's theorem we have
\begin{equation}
\mapleinline{inert}{2d}{St := Sum(GAMMA(k)/((k+a-1)*GAMMA(k+a)), k = 1 .. infinity) = (1/2)*(Psi(1, (1/2)*a)-Psi(1, (1/2)*a+1/2))/GAMMA(a)}{\[\displaystyle \sum _{k=1}^{\infty }{\frac {\Gamma \left( k \right) }{ \left( k+a-1 \right) \Gamma \left( k+a \right) }}={\frac {\Psi^{\prime} \left(a/2 \right) -\Psi^{\prime} \left( a/2+1/2 \right) 
\mbox{}}{2\,\Gamma \left( a \right) }}\]},
\label{Lemma2b}
\end{equation}
and, after differentiating with respect to $a$, we obtain

\begin{align} \label{Lemma2}
&\displaystyle \sum _{k=1}^{\infty }{\frac {\Gamma \left( k \right) \Psi \left( k+a \right) }{ \left( k+a-1 \right) \Gamma \left( k+a \right) }}=-\sum _{k=1}^{\infty }{\frac {\Gamma \left( k \right) }{ \left( k+a-1 \right) ^{2}\Gamma \left( k+a \right) }}
\mbox{}  \\ \nonumber
&+{\frac { \left( \Psi^{\prime} \left(a/2 \right) -\Psi^{\prime} \left( a/2+1/2 \right)  \right) \Psi \left( a \right)/2 -1/4\,\Psi^{(2)} \left( a/2 \right) 
\mbox{}+1/4\,\Psi^{(2)} \left(a/2+1/2 \right) }{\Gamma \left( a \right) }}
\end{align}

\item{\bf Lemma 3}

By adding and subtracting a term corresponding to $k=-1$, from its series representation, we have
\begin{equation}
\mapleinline{inert}{2d}{hypergeom([1, 2, 2-a, a+1], [3, 3, 3], 1) = 8*Pi*Hypergeom([1, 1, a, 1-a], [2, 2, 2], 1)/((1-a)*GAMMA(1-a)*GAMMA(a)*a*sin(Pi*a))-8/((1-a)*a)}{\[\displaystyle {\mbox{$_4$F$_3$}(1,2,2-a,a+1;\,3,3,3;\,1)}={\frac {8\,\pi\,{\mbox{$_4$F$_3$}}\left( 1,1,a,1-a;2,2,2;1 \right) }{ \left( 1-a \right) \Gamma \left( 1-a \right) 
\mbox{}\Gamma \left( a \right) a\sin \left( \pi\,a \right) }}-{\frac {8}{ \left( 1-a \right) a}}\]}
\label{Lemma3a}
\end{equation}
so that, after the invocation of \eqref{Y30} 

\begin{align} \label{Lemma3b}
&\displaystyle {\mbox{$_4$F$_3$}(1,2,2-a,a+1;\,3,3,3;\,1)}={\frac {\sin \left( \pi\,a \right) }{{a}^{2} \left( a-1 \right) ^{2}
\mbox{}\pi} \left( 16\,\Psi^{\prime} \left(a \right) -8\,\Psi^{\prime} \left( a/2 \right) +8\,{\frac {{a}^{2}-a+1}{{a}^{2} \left( a-1 \right) ^{2}}} \right) }\\ \nonumber
&-16\,{\frac {\Psi \left( a+1 \right) }{{a}^{2} \left( a-1 \right) ^{2}}}
-16\,{\frac {\gamma}{{a}^{2} \left( a-1 \right) ^{2}}}+8\,{\frac {{a}^{4}-2\,{a}^{3}+{a}^{2}+2\,a-1}{{a}^{3} \left( a-1 \right) ^{3}}}\,.
\end{align}

\item{\bf Lemma 4}

From \cite[Entry 26]{Milgram447} after evaluating some limits,
\begin{align} \label{Lemma4}
&\displaystyle {\mbox{$_3$F$_2$}(1,1,n;\,n+1,n+1;\,1)}=\\ \nonumber
&- \left( \sum _{k=0}^{n-3}{\frac { \left( -1 \right) ^{k}\Psi \left( 1+k \right) }{\Gamma \left( n-k-1 \right)  \left( n-k-1 \right) ^{2}\Gamma \left( 1+k \right) 
\mbox{}}}+\sum _{k=0}^{n-3}{\frac { \left( -1 \right) ^{k}}{\Gamma \left( n-k-1 \right)  \left( n-k-1 \right) ^{3}\Gamma \left( 1+k \right) }} \right) {n}^{2}\Gamma \left( n \right) \\ \nonumber 
&+ \left( -1 \right) ^{n} \left( 1/2- \left( \Psi \left( n \right) +1 \right) {n}^{3}+ {n}^{2}\,\bigg{(} {\gamma}^{2}/2\,+2\,\gamma\,\Psi \left( n \right) +{\pi}^{2}/12\,+3/2\,\Psi \left( n \right) ^{2}+\Psi \left( n \right) \right. \\ \nonumber
& \left.  \hspace{2cm}  -2\,\Psi^{\prime} \left(n \right) +1/2\,\Psi^{\prime} \left(n+1  \right)  
+2 \bigg{)} \right) 
\end{align}

\item{\bf Lemma 5}

From \cite[Entry 13]{Milgram447} - contiguous to Whipples theorem, and a generalization of \eqref{Lemma4}
\begin{equation}
\mapleinline{inert}{2d}{Ht := hypergeom([1, 1, a], [a+1, a+1], 1) = (1/4)*a^2*Psi(1, (1/2)*a+1)+(1/4)*a^2*Psi(1, (1/2)*a)-(1/2)*a^2*Psi(1, (1/2)*a+1/2)+1}{\[\displaystyle {\mbox{$_3$F$_2$}(1,1,a;\,a+1,a+1;\,1)}=\frac{{a}^{2}}{4}\Psi^{\prime} \left(a/2+1 \right) +\frac{{a}^{2}}{4}\Psi^{\prime} \left(a/2 \right) -\frac{{a}^{2}}{2}\Psi^{\prime} \left(a/2+1/2 \right) 
\mbox{}+1\]}
\label{ht}
\end{equation}
so that, after differentiating we obtain

\begin{align}\label{Lemma5}
&\displaystyle \sum _{k=1}^{\infty }{\frac {\Gamma \left( k \right) \Psi \left( k+a \right) }{ \left( k+a-1 \right) \Gamma \left( k+a \right) }}=-\sum _{k=1}^{\infty }{\frac {\Gamma \left( k \right) }{ \left( k+a-1 \right) ^{3}\Gamma \left( k+a-1 \right) }}
\mbox{}\\ \nonumber
&+{\frac {2\,\Psi^{\prime} \left(a/2 \right) \Psi \left( a \right) -2\,\Psi^{\prime} \left( a/2+1/2 \right) \Psi \left( a \right) -\Psi^{(2)} \left( a/2 \right) +\Psi^{(2)} \left(a/2+1/2 \right) 
}{4\,\Gamma \left( a \right) }}
\end{align}
For a comparison of \eqref{Lemma4} and \eqref{ht}, see \eqref{Newsum2}.

\item{\bf Lemma 6}

Trivially,
\begin{equation} 
{\mbox{$_3$F$_2$}(1,1,2-a;\,2,2;\,1)}=\frac{\Psi(a)+\gamma}{\Gamma(a)}
\end{equation}
leading to

\begin{equation}
\mapleinline{inert}{2d}{Sum((-1)^k*Psi(a-k)/(k^2*GAMMA(k)*GAMMA(a-k)), k = 1 .. infinity) = Psi(1, a)/GAMMA(a)-Psi(a)^2/GAMMA(a)-Psi(a)*gamma/GAMMA(a)}{\[\displaystyle \sum _{k=1}^{\infty }{\frac { \left( -1 \right) ^{k}\Psi \left( a-k \right) }{{k}^{2}\Gamma \left( k \right) \Gamma \left( a-k \right) }}={\frac {\Psi^{\prime} \left( a \right) }{\Gamma \left( a \right) }}-{\frac {\Psi \left( a \right)^{2}}{\Gamma \left( a \right) }}
\mbox{}-\gamma{\frac {\Psi \left( a \right)}{\Gamma \left( a \right) }}\]}
\label{Lemma6}
\end{equation}

\end{itemize}
\end{appendices}

\end{flushleft}
\end{document}